\documentclass[11pt,reqno,psamsfonts]{amsart}
\pdfoutput=1

\usepackage{amssymb}
\usepackage{amsthm}
\usepackage{amsmath}
\usepackage{latexsym}
\usepackage[T1]{fontenc}
\usepackage[utf8]{inputenc}
\usepackage[russian, french, 
english]{babel}

\usepackage{graphicx}
\usepackage{wrapfig}
\usepackage{mathtools}
\usepackage{amsbsy}
\usepackage[inline]{enumitem}
\usepackage{mathrsfs}
\usepackage{array}
\usepackage{multicol}
\usepackage{stmaryrd}
\usepackage{cancel}
\usepackage{lmodern}
\usepackage{mathabx}
\usepackage{upgreek}
\usepackage{titlesec}
\usepackage{titletoc}
\usepackage[spacing=true,kerning=true,babel=true,tracking=true]{microtype}
\usepackage[shortcuts]{extdash}
\usepackage[foot]{amsaddr}
\usepackage[left=1in,right=1in,top=1in,bottom=1in,bindingoffset=0cm]{geometry}
\usepackage{bm}
\usepackage{centernot}
\usepackage{mdframed}
\usepackage[hidelinks]{hyperref}
\usepackage{xspace}
\usepackage[most]{tcolorbox}
\usepackage{framed}
\usepackage[labelfont=bf]{caption}
\usepackage{tikz}
\usetikzlibrary{math}
\usetikzlibrary{shapes,snakes}
\usetikzlibrary{arrows.meta}
\usetikzlibrary{decorations.pathmorphing}
\usetikzlibrary{patterns}
\usetikzlibrary{positioning,calc,intersections}
\usepackage{tikzsymbols}
\usepackage{float}

\DeclareFontFamily{U}{skulls}{}
\DeclareFontShape{U}{skulls}{m}{n}{ <-> skull }{}

\usepackage[
    sortcites,
    backend=biber, style=alphabetic, sorting=nyt, maxnames=100,backref=true]{biblatex}
\newtheoremstyle{bfnote}%
{}{}%
{\slshape}{}%
{\bfseries}{\bfseries.}%
{ }%
{\thmname{#1}\thmnumber{ #2}\thmnote{ \ep{\normalfont{}#3}}}

\theoremstyle{bfnote}
\newtheorem{theo}{Theorem}[section]

\newtheorem*{theo*}{Theorem}
\newtheorem{prop}[theo]{Proposition}
\newtheorem{lemma}[theo]{Lemma}
\newtheorem{claim}[theo]{Claim}

\newtheorem{corl}[theo]{Corollary}

\newtheorem*{corl*}{Corollary}

\theoremstyle{definition}
\newtheorem{defn}[theo]{Definition}
\newtheorem*{defn*}{Definition}
\newtheorem{exmp}[theo]{Example}

\newtheorem{remk}[theo]{Remark}

\newtheorem*{remks*}{Remarks}
\newtheorem*{exmp*}{Example}

\theoremstyle{remark}
\newtheorem*{ques*}{Question}
\newtheorem*{remk*}{Remark}

\newcommand*{\myproofname}{Proof}
\newenvironment{claimproof}[1][\myproofname]{\begin{proof}[#1]}{\end{proof}}

\makeatletter
\newcommand{\neutralize}[1]{\expandafter\let\csname c@#1\endcsname\count@}
\makeatother

\newenvironment{theocopy}[1]
{\renewcommand{\thetheo}{\ref{#1}}%
	\neutralize{theo}\phantomsection
	\begin{theo}}
	{\end{theo}}

\newenvironment{lemmacopy}[1]
{\renewcommand{\thetheo}{\ref{#1}}%
	\neutralize{theo}\phantomsection
	\begin{lemma}}
	{\end{lemma}}

\newcommand{\0}{\varnothing}
\newcommand{\set}[1]{\{#1\}}

\renewcommand{\P}{\mathbb{P}}

\renewcommand{\epsilon}{\varepsilon}

\renewcommand{\phi}{\varphi}
\renewcommand{\theta}{\vartheta}
\renewcommand{\leq}{\leqslant}
\renewcommand{\geq}{\geqslant}
\newcommand{\defeq}{\coloneqq}

\newcommand{\bemph}[1]{{\normalfont#1}}
\newcommand{\ep}[1]{\bemph{(}#1\bemph{)}}

\newcommand{\emphd}[1]{{\fontseries{b}\selectfont\textsf{#1}}}
\newcommand{\acts}{\mathrel{\reflectbox{$\righttoleftarrow$}}}
\newcommand{\G}{\Gamma}
\newcommand{\pto}{\dashrightarrow}
\newcommand{\Free}{\mathsf{Free}}

\newcommand{\symdif}{\mathbin{\triangle}}
\newcommand{\rest}[2]{{{#1}\vert_{#2}}}

\newcommand{\fs}[2]{{[{#1}]^{#2}}}
\newcommand{\finset}[1]{\fs{{#1}}{{<\infty}}}
\newcommand{\N}{\mathbb{N}}
\newcommand{\Z}{\mathbb{Z}}

\newcommand{\dom}{\mathsf{dom}}
\newcommand{\asi}{\mathsf{asi}}
\newcommand{\casi}{\asi_{\mathsf{c}}}

\newcommand{\Cay}{\mathsf{Cay}}
\newcommand{\Sch}{\mathsf{Sch}}
\newcommand{\Sub}{\mathsf{Sub}}
\newcommand{\In}{\mathsf{In}}
\newcommand{\Out}{\mathsf{Out}}
\newcommand{\iden}{\mathbf{1}}
\newcommand{\sr}{amply syndetic\xspace}
\newcommand{\Sr}{Amply syndetic\xspace}
\newcommand{\Sep}{\mathsf{Sep}}
\newcommand{\Synd}{\mathsf{Synd}}
\newcommand{\Min}{\mathsf{Min}}
\newcommand{\Vis}{\mathsf{Vis}}

\numberwithin{equation}{section}

\newenvironment{scproof}[1][]{\begin{proof}[\textsc{\upshape{Proof}}#1]}{\end{proof}}

\titleformat{\section}[block]{\large\bfseries\sffamily}{\thesection.}{1ex}{}
\titleformat{\subsection}[block]{\bfseries\sffamily}{\thesubsection.}{1ex}{}
\titleformat{\subsubsection}[runin]{\bfseries}{\bfseries\upshape\thesubsubsection.}{1ex}{}[.]

\titlespacing*{\section}{0pt}{*3}{*1}
\titlespacing*{\subsection}{0pt}{*3}{*1}
\titlespacing*{\subsubsection}{0pt}{*1.5}{*1}

 \titlecontents{section}
 [1.5em] %
 {\smallskip}
 {\bfseries\thecontentslabel\hspace{1.02em}}
{\bfseries}
 {\,\,\titlerule*[0.77pc]{}\bfseries\contentspage}
\titlecontents{subsection}
 [4em] %
 {\smallskip}
 {\thecontentslabel\hspace{1.02em}}
{\hspace*{2.32em}}
 {\,\,\titlerule*[0.77pc]{.}\contentspage}

\renewbibmacro{in:}{}

\renewbibmacro*{volume+number+eid}{%
	\printfield{volume}%
	\setunit*{\addnbspace}
	\printfield{number}%
	\setunit{\addcomma\space}%
	\printfield{eid}}

\DeclareFieldFormat[article]{volume}{\textbf{#1}\space}
\DeclareFieldFormat[article]{number}{\mkbibparens{#1}}

\DeclareFieldFormat{journaltitle}{#1,}
\DeclareFieldFormat[thesis]{title}{\mkbibemph{#1}\addperiod}
\DeclareFieldFormat[article, unpublished, thesis]{title}{\mkbibemph{#1},}
\DeclareFieldFormat[book]{title}{\mkbibemph{#1}\addperiod}
\DeclareFieldFormat[unpublished]{howpublished}{#1, }

\DeclareFieldFormat{pages}{#1}

\DeclareFieldFormat[article]{series}{Ser.~#1\addcomma}

\setlist{topsep=3pt,itemsep=3pt}

\title{\sffamily Flows with minimal subdynamics}
\date{}

\author{Anton~Bernshteyn}
\address{\normalfont (AB) Department of Mathematics, University of California, Los Angeles, CA, USA}
\email{bernshteyn@math.ucla.edu}

\author{Joshua~Frisch}
\address{\normalfont (JF) 
Department of Mathematics, University of California, San Diego, CA, USA}
\email{jfrisch@ucsd.edu}

\thanks{AB's research is partially supported by the NSF CAREER grant DMS-2528522 and the Sloan Research Fellowship (2025).
JF's research is partially supported by the NSF grant DMS-2348981.}

\begin{document}


\maketitle

\begin{abstract}
    Let $\G$ be a countably infinite discrete group. A $\G$-flow $X$ (i.e., a nonempty compact Hausdorff space equipped with a continuous action of $\G$) is called {$S$-minimal} for a subset $S \subseteq \G$
    if the partial orbit $S \cdot x$ is dense for every point
    $x \in X$. We show that for any countable family $(S_n)_{n \in \N}$ of infinite subsets of $\G$, there exists a free $\G$-flow $X$ that is $S_n$-minimal for all $n \in \N$; additionally, $X$ can be taken to be a subflow of $2^\G$. This vastly generalizes a result of Frisch, Seward, and Zucker,
    in which each $S_n$ is required to be a normal subgroup of $\G$. As a corollary, we show that for a given Polish $\G$-flow $X$, there exists a free $\G$-flow $Y$ disjoint from $X$ in the sense of Furstenberg if and only if $X$ has no wandering points. This completes a line of inquiry started by Glasner, Tsankov, Weiss, and Zucker. 
    As another application, we strengthen some of the results of Gao, Jackson, Krohne, and Seward on the structure of Borel complete sections. For example, we show that if $B$ is a Borel complete section in the free part of $2^\G$, then every union of sufficiently many shifts of $B$ contains an orbit (previously, this was only known for open sets $B$). Although our main results are purely dynamical, their proofs rely on recently developed machinery from descriptive set-theoretic combinatorics, namely the 
    asymptotic separation index introduced by Conley, Jackson, Marks, Seward, and Tucker-Drob and its links to the Lov\'asz Local Lemma.
\end{abstract}

\section{Introduction}\label{sec:intro}

    \subsection{Existence of flows and subshifts with minimal subdynamics}

    Throughout, $\G$ is an infinite
    countable discrete group.
    A \emphd{$\G$-flow} is a nonempty compact Hausdorff space $X$ equipped with a continuous action $\G \acts X$. A \emphd{subflow} of a $\G$-flow $X$ is a nonempty closed $\G$-invariant subset $Y \subseteq X$. A $\G$-flow $X$ is \emphd{minimal} if it has no subflows other than $X$ itself; equivalently, $X$ is minimal if the orbit $\G \cdot x$ of every point $x \in X$ is dense in $X$. A routine argument using Zorn's lemma shows that every $\G$-flow has a minimal subflow. That being said, finding explicit examples of minimal $\G$-flows with desirable properties is often a challenging task, leading to many intriguing open questions. In this paper we are interested in the following class of problems, which Seward, Zucker, and the second named author \cite{FSZ} called \emph{minimal subdynamics}:

    \begin{tcolorbox}[width=\linewidth, sharp corners=all, colback=white!95!black, boxrule=0pt,colframe=white] 
        Given a family $(\Delta_i)_{i \in I}$ of infinite subgroups of $\G$, when does there exist a free $\G$-flow $X$ such that the induced action $\Delta_i \acts X$ makes it a minimal $\Delta_i$-flow for all $i \in I$?
    \end{tcolorbox}

    Here we say that a group action $\G \acts X$ is \emphd{free} if the stabilizer of every point $x \in X$ is trivial, i.e., if $\gamma \cdot x \neq x$ for all $x \in X$ and $\gamma \in \G \setminus \set{\iden}$, where $\iden$ is the identity element of $\G$.

    The study of dynamical systems via their ``subdynamics,'' i.e., by restricting the action to a subgroup, 
    has long been a theme 
    in topological dynamics and ergodic theory. A classical example is the notion of a \emph{totally minimal} homeomorphism, i.e., a homeomorphism $T \colon X \to X$ such that for all nonzero $n \in \Z$, the action $\langle T^n \rangle \acts X$ of the group generated by $T^n$ is minimal. Totally minimal homeomorphisms were introduced by Hedlund in \cite{HedlundSturmian} (under the name ``powerfully minimal'').  
    This definition is extended to $\G$-flows for general $\G$ by saying that a $\G$-flow is \emph{totally minimal} if it is minimal as a $\Delta$-flow for all subgroups $\Delta \leq \G$ of finite index \cite[Defn.~2.27]{GottHed}. \emph{Totally ergodic} is an analogous concept in ergodic theory \cite[Defn.~10.3]{RSErgodic}. Various aspects of subdynamics for $\Z^d$-flows, such as the notion of \emph{expansive subdynamics} introduced 
    by Boyle and Lind in \cite{Expansive1}, have been extensively studied; see, e.g., \cite{Expansive2,Expansive3,Expansive4,Expansive5,Expansive6,Expansive7,Expansive8,Expansive9,Rosenthal}  for a small sample of the literature on the subject.
    
    Generalizing an earlier result of Zucker \cite{Zucker}, Seward, Zucker, and the second named author gave the following partial answer to the minimal subdynamics question:

    \begin{theo}[{JF--Seward--Zucker \cite[Thm.~0.3]{FSZ}}]\label{theo:FSZ}
        If $(\Delta_n)_{n \in \N}$ is a countable family of infinite {normal} subgroups of $\G$, then
        there is a free $\G$-flow that is minimal as a $\Delta_n$-flow for all $n \in \N$.
    \end{theo}

    The normality assumption on the subgroups $\Delta_n$ is quite restrictive, but it appears difficult to remove it with prior techniques.
    Indeed, even the following very modest special case of the general problem remained open until now:

    \begin{quote}
        {Let $\mathbb{F}_2$ be the free group of rank $2$. Is it true that for every non-identity element $\gamma$ of $\mathbb{F}_2$, 
        there exists a free $\mathbb{F}_2$-flow that is minimal as a $\langle \gamma \rangle$-flow?}
    \end{quote}

    In this paper we introduce a novel approach that combines tools from topological dynamics, descriptive set theory, and combinatorics and allows us to completely remove the normality assumption in Theorem~\ref{theo:FSZ}. Not only that,  
    we do not even need to assume that each $\Delta_n$ is a \emph{subgroup} of $\G$. 

    \begin{defn}[$S$-minimal $\G$-flows]\label{defn:S-minimal}
        We say that a $\G$-flow $X$ is \emphd{$S$-minimal} for a {subset} $S \subseteq \G$ if the partial orbit $S \cdot x$ is dense in $X$ for every point $x \in X$.
    \end{defn}


    \begin{tcolorbox}
    \begin{theo}\label{theo:main}
        If $(S_n)_{n \in \N}$ is a countable family of infinite subsets of $\G$, then there exists a free $\G$-flow 
        that is $S_n$-minimal for all $n \in \N$.
        %
        %
        %
        %
    \end{theo}
    \end{tcolorbox}



    
    The following statements are immediate consequences of Theorem~\ref{theo:main}:

    \begin{tcolorbox}
    \begin{corl}\label{corl:fg_subgroups}
        There exists a free $\G$-flow 
        that is $\Delta$-minimal with respect to
        every non-locally finite subgroup $\Delta \leq \G$.
    \end{corl}
    \end{tcolorbox}
    \begin{scproof}
        Note that if $S \subseteq S' \subseteq \G$, then every $S$-minimal $\G$-flow is also $S'$-minimal. Hence, it suffices to find a free $\G$-flow 
        that is $\Delta$-minimal for every infinite {finitely generated} subgroup $\Delta \leq \G$. Such a $\G$-flow 
        exists by Theorem~\ref{theo:main} as $\G$ has only countably many finitely generated subgroups. 
    \end{scproof}



    \begin{tcolorbox}
    \begin{corl}\label{corl:no_loc_fin}
        If $\G$ has at most countably many infinite locally finite subgroups, then there exists a free $\G$-flow 
        that is $\Delta$-minimal with respect to every infinite subgroup $\Delta \leq \G$.
    \end{corl}
    \end{tcolorbox}
    \begin{scproof}
        As in the proof of Corollary~\ref{corl:fg_subgroups}, we just need to find a free $\G$-flow 
        that is $\Delta$-minimal for every infinite subgroup $\Delta \leq \G$ that is either finitely generated or locally finite. 
    \end{scproof}

    In particular, there exists a free $\mathbb{F}_2$-flow that is $\langle \gamma \rangle$-minimal for all non-identity elements $\gamma \in \mathbb{F}_2$ simultaneously. In general, Corollaries~\ref{corl:fg_subgroups} and \ref{corl:no_loc_fin} are new 
    for most non-Abelian groups $\G$. 
    Notice that a nontrivial $\G$-flow cannot be $S$-minimal for all infinite subsets $S \subseteq \G$ at once; in that sense, the restriction to subgroups in Corollary~\ref{corl:no_loc_fin} is essential.

    \begin{remk}
        The conclusion of Corollary~\ref{corl:no_loc_fin} may fail if $\G$ has uncountably many infinite locally finite subgroups. 
    For example,
    Seward, Zucker, and the second named author showed that  it fails for 
    $\G = \bigoplus_{n \in \N} (\Z/2\Z)$ 
    \cite[Thm.~2.5]{FSZ}. 
    Exactly characterizing the class of groups $\G$ that satisfy 
    the conclusion of Corollary~\ref{corl:no_loc_fin} 
    remains an interesting open problem.
    \end{remk}



    We can generalize Theorem~\ref{theo:main} further using the following definitions:

    \begin{defn}[Syndetic sets]
        Let $\G \acts X$ be an action of $\G$. A set $U \subseteq X$ is \emphd{$F$-syndetic} for a finite subset $F \subset \G$ if $F^{-1} \cdot U = X$, i.e., if for each $x \in X$, there is $\sigma \in F$ such that $\sigma \cdot x \in U$.
    \end{defn}

    \begin{defn}[$\mathcal{F}$-minimal $\G$-flows]\label{defn:F-minimal}
        Let $\mathcal{F}$ be a family of finite subsets of $\G$. We say that a $\G$-flow $X$ is \emphd{$\mathcal{F}$-minimal} if for every nonempty open set $U \subseteq X$, there is $F \in \mathcal{F}$ such that $U$ is $F$-syndetic. 
    \end{defn}

    Definition~\ref{defn:F-minimal} extends Definition~\ref{defn:S-minimal} thanks to the following well-known observation:

    \begin{prop}\label{prop:synd1}
        A $\G$-flow $X$ is $S$-minimal for a subset $S \subseteq \G$ if and only if it is $\finset{S}$-minimal, where $\finset{S}$ is the family of all finite subsets of $S$. 
    \end{prop}
    \begin{scproof}
        Clearly, if $X$ is $\finset{S}$-minimal, the partial orbit
        $S \cdot x$ is dense for all $x \in X$. Conversely, suppose $X$ is $S$-minimal and let $U \subseteq X$ be nonempty and open. Then we have  
        $S^{-1} \cdot U = X$, and the compactness of $X$ yields a finite set $F \in \finset{S}$ such that $F^{-1}\cdot U = X$. 
    \end{scproof}

    A family of finite sets $\mathcal{F}$ is called \emphd{unbounded} if $\sup_{F \in \mathcal{F}} |F| = \infty$. Theorem~\ref{theo:main} 
    remains valid if we replace the sets $S_n$ and the corresponding families $\finset{S_n}$ by {arbitrary} unbounded families $\mathcal{F}_n$: 

    \begin{tcolorbox}
    \begin{theo}\label{theo:main_prime}
        If $(\mathcal{F}_n)_{n \in \N}$ is a sequence of unbounded families of finite subsets of $\G$, then there exists a free $\G$-flow 
        that is $\mathcal{F}_n$-minimal for all $n \in \N$.
    \end{theo}
    \end{tcolorbox}


    From now on, we shall treat Theorem~\ref{theo:main_prime} as our main result (with Theorem~\ref{theo:main} as a special case), and most of the paper will be devoted to its proof. Although the concept of $\mathcal{F}$-minimality may appear somewhat artificial, it provides a natural framework for our proof techniques. Additionally, for some of the 
    applications we present here, 
    the general setting of $\mathcal{F}$-minimality is essential. 


    %
    %
    
    Our proof of Theorem~\ref{theo:main_prime} naturally yields a $\G$-flow whose underlying space is Polish---indeed, it is homeomorphic to the Cantor space. With a little extra work, we are able to give it a particularly nice special form, namely that of a \emph{subshift}. 
    We use the notation 
    $\N \defeq \set{0,1,2,\ldots}$ and $\N^+\defeq\set{1,2,3,\ldots}$. 
    Each natural number
    $k \in \N$ is identified with the $k$-element set $\set{i \in \N \,:\, i < k}$ and given the discrete topology. For $k \in \N^+$, the \emphd{shift action} $\G \acts k^\G$ 
    of $\G$ on the product space $k^\G$ is defined by
    \[
        (\gamma \cdot x)(\delta) \,\defeq\, x(\delta \gamma) \quad \text{for all $x \colon \G \to k$ and $\gamma$, $\delta \in \G$}.
    \]
    This makes $k^\G$ a 
    $\G$-flow, called a \emphd{Bernoulli shift}, or simply a \emphd{shift}. Subflows of $k^\G$ are called \emphd{subshifts} and are important examples of dynamical systems 
    studied in {symbolic dynamics} \cite{subshifts1,subshifts2}. 

    \begin{tcolorbox}
    \begin{corl}\label{corl:subshift}\label{corl:subshift_prime}
        If $(\mathcal{F}_n)_{n \in \N}$ is a sequence of unbounded families of finite subsets of $\G$, then there exists a free subflow of $2^\G$
        that is $\mathcal{F}_n$-minimal for all $n \in \N$.
        %
        %
        %
        %
    \end{corl}
    \end{tcolorbox}

    In particular, for any countable family $(S_n)_{n \in \N}$ of infinite subsets of $\G$, there exists a free subflow of $2^\G$ 
    that is $S_n$-minimal for all $n \in \N$.    
    To derive Corollary~\ref{corl:subshift} from Theorem~\ref{theo:main_prime}, we use a result of Seward and Tucker-Drob \cite{ST-D}, which yields a free subshift $Y \subset 2^\G$ with a certain ``universality'' property, and then find the desired $\G$-flow among the subflows of $Y$. 
    See \S\ref{sec:subshift} for details.
    
    We should note that the existence of a free subshift for an arbitrary countable group $\G$ is already a highly nontrivial fact. It was established in full generality by Gao, Jackson, and Seward \cite{GJS1,GJS2} after partial results in various special cases due to Dranishnikov and Schroeder \cite{DS} and Glasner and Uspenskij \cite{GU} (the $\G = \Z$ case appears implicitly as far back as the work of Thue \cite{Thue}, Morse \cite{Morse}, and Morse--Hedlund \cite{MorseHedlund}). 
    Aubrun, Barbieri, and Thomass\'{e} \cite{ABT} subsequently found a simpler proof for general $\G$ using the probabilistic method. Since then, several works have appeared that construct free subshifts with various additional properties \cite{ST-D,LLLDyn1,LLLDyn4}, and Corollary~\ref{corl:subshift} continues this trend.

    \subsection{Applications to constructing disjoint flows}

    It turns out that our results have close ties to the theory of {disjointness} for $\G$-flows, introduced by Furstenberg in his seminal paper \cite{Furstenberg}. 

    \begin{defn}[Joinings and disjointness]
        Let $X$ and $Y$ be $\G$-flows. We view $X \times Y$ as a $\G$-flow with the diagonal action of $\G$. A \emphd{joining} of $X$ and $Y$ is a subflow $Z \subseteq X \times Y$ that projects onto $X$ and $Y$. The $\G$-flows $X$, $Y$ are \emphd{disjoint}, in symbols $X \perp Y$, if their only joining is $X \times Y$ itself.
    \end{defn}

    An easy observation is that if $X \perp Y$, then at least one of $X$, $Y$ is minimal. Also, if both $X$ and $Y$ are minimal, then $X \perp Y$ if and only if $X \times Y$ is also minimal.

    Every $\G$-flow is disjoint from the trivial action of $\G$ on a single-point space, but this is hardly an interesting instance of the disjointness relation. Seeking more illuminating examples, we ask: 

    \begin{tcolorbox}[width=\linewidth, sharp corners=all, colback=white!95!black, boxrule=0pt,colframe=white] 
        For which $\G$-flows $X$ does there exist some free $\G$-flow $Y$ such that $X \perp Y$?
    \end{tcolorbox}

    We shall specifically focus on {Polish} (equivalently, metrizable) $\G$-flows $X$. Some important results related to this question were obtained by Glasner, Tsankov, Weiss, and Zucker in their breakthrough paper \cite{BernoulliDisjointness}. Generalizing a theorem of Furstenberg \cite{Furstenberg}, they showed that $2^\G$ and, more generally, all so-called \emph{strongly irreducible subshifts} are disjoint from every minimal $\G$-flow \cite[Thm.~6.2]{BernoulliDisjointness}. They also showed that if a Polish $\G$-flow $X$ is minimal, then it is disjoint from some free minimal Polish $\G$-flow \cite[Thm.~1.2(i)]{BernoulliDisjointness}. The assumption that $X$ is Polish cannot be removed here, as the \emph{universal minimal flow} $\mathcal{M}(\G)$ is a non-metrizable counterexample. 
    
    Using our results on minimal subdynamics, we give a complete answer to the above question (for Polish $\G$-flows). We need the following standard definition:

    \begin{defn}[Wandering sets and points]
        A subset $U \subseteq X$ of a $\G$-flow $X$ is \emphd{wandering} if $U \cap (\gamma \cdot U) = \0$ for all but finitely many group elements $\gamma \in \G$, and a point $x \in X$ is \emphd{wandering} if it has a wandering neighborhood.
    \end{defn}

    By definition, the set of all wandering points in a $\G$-flow $X$ is open and $\G$-invariant. It is clear that if $x \in X$ is a wandering point, then the stabilizer of $x$ is finite. It is also not hard to see that if the stabilizer of a wandering point $x$ is trivial (for instance, if $\G$ is torsion-free), then $x$ has a neighborhood $U$ such that 
    the sets $(\gamma \cdot U)_{\gamma \in \G}$ are pairwise disjoint. 

    We show that the existence of wandering points is the only obstruction to being disjoint from some free $\G$-flow:

    \begin{tcolorbox}
    \begin{theo}\label{theo:disjoint_one}
        The following statements are equivalent for a Polish $\G$-flow $X$:
        \begin{enumerate}[label=\ep{\normalfont\arabic*}]
            \item $X$ is disjoint from some infinite $\G$-flow,
            \item $X$ is disjoint from some free minimal subflow of $2^\G$,
            \item $X$ has no wandering points.
        \end{enumerate}
    \end{theo}
    \end{tcolorbox}

    Furthermore, we can extend this result to countable families of Polish $\G$-flows:

    \begin{tcolorbox}
    \begin{theo}\label{theo:disjoint_many}
        If $(X_n)_{n \in \N}$ is a countable family of Polish $\G$-flows with no wandering points, then there exists a free minimal subflow of $2^\G$ that is disjoint from $X_n$ for all $n \in \N$.
    \end{theo}
    \end{tcolorbox}

    These results are established in \S\ref{sec:disjoint} by reducing the disjointness relation to an instance of the minimal subdynamics problem (see Proposition~\ref{prop:characterize}). 
    It is possible to deduce Theorem~\ref{theo:disjoint_many} from Theorem~\ref{theo:disjoint_one} by letting $X$ be the one-point compactification of 
    $\bigsqcup_{n \in \N} X_n$, but we give a direct proof instead. Minimal $\G$-flows have no wandering points (see, e.g., \cite[Lem.~2.1]{LLLDyn3}), so Theorem~\ref{theo:disjoint_one} includes the aforementioned result \cite[Thm.~1.2(i)]{BernoulliDisjointness} of Glasner \emph{et al.} as a special case.  

    One of the reasons these facts have not been discovered earlier is their 
    tight connection to the minimal subdynamics problem. 
    Indeed, consider the following example. 
    Let $\Delta \leq \G$ be a subgroup. Then $\G$ naturally acts on the set $\G/\Delta$ of the left cosets of $\Delta$. We endow $\G/\Delta$ with the discrete topology and let $(\G/\Delta)^*$ be the one-point extension of $\G/\Delta$, where the action of $\G$ is extended to the point at infinity by making it fixed. (Note that if the index of $\Delta$ in $\G$ is finite, then the point at infinity is isolated.) 
    %
    The space $(\G/\Delta)^*$ is a (countable) Polish $\G$-flow, and we have the following:

    \begin{prop}\label{prop:cosets}
        %
        The following statements are equivalent for a $\G$-flow $X$ and a subgroup $\Delta \leq \G$:
        \begin{enumerate}[label=\ep{\normalfont\arabic*}]
            \item $X \perp (\G/\Delta)^*$,
            \item $X$ is $\Delta$-minimal.
        \end{enumerate}
    \end{prop}

    See \S\ref{sec:disjoint} for the proof. If $\Delta$ is infinite, then there are no wandering points in $(\G/\Delta)^*$ (because the stabilizer of every point is infinite). Therefore, Theorem~\ref{theo:disjoint_many} implies the subgroup case of Theorem~\ref{theo:main}, i.e., it shows that for any countable family $(\Delta_n)_{n \in \N}$ of infinite subgroups of $\G$, there exists a free $\G$-flow $X$ that is $\Delta_n$-minimal for all $n \in \N$. In particular, Theorem~\ref{theo:disjoint_many} is already sufficient to derive Corollaries~\ref{corl:fg_subgroups} and \ref{corl:no_loc_fin}.

    \subsection{Applications to Borel complete sections}\label{subsec:sections}

    As mentioned previously, the existence of free subshifts for arbitrary $\G$ was established by Gao, Jackson, and Seward \cite{GJS1,GJS2}. Their work 
%
 was motivated by the realization that free subshifts offer deep insight into the topological and descriptive set-theoretic properties of 
 the shift action $\G \acts 2^\G$ 
 by putting significant constraints on the geometry of complete sections in $\Free(2^\G)$ (see below for the definition). This idea was applied with great success in subsequent contributions by
    %
    Gao, Jackson, Krohne, and Seward \cite{GJKS_Borel,GJKS_continuous}. 
    %
    %
    It turns out that some of the 
    main results of 
    \cite{GJKS_Borel}
    can be significantly improved with the help of free subshifts with minimal subdynamics, 
    specifically using the notion of $\mathcal{F}$-minimality (Definition~\ref{defn:F-minimal}).
    
    Let $\Free(2^\G)$ denote the \emphd{free part} of $2^\G$, i.e., the set of all points $x \in 2^\G$ 
    with trivial stabilizers.
    Then $\Free(2^\G)$ is a shift-invariant dense $G_\delta$ subset of $2^\G$ 
    (see, e.g.,
    \cite[Lem.~2.3]{BernoulliDisjointness} for a proof of density), and it is the largest subspace of $2^\G$ on which $\G$ acts freely. 
    A \emphd{complete section} in $\Free(2^\G)$ is a subset $B \subseteq \Free(2^\G)$ that meets every $\G$-orbit, i.e.,
    such that
    $\G \cdot B = \Free(2^\G)$. 
    The following is a typical fact that can be proved using Gao \emph{et al.}'s methods:

    \begin{theo}[{Gao--Jackson--Krohne--Seward \cite[Thm.~1.1]{GJKS_Borel}}]\label{theo:GJKS1}
        Suppose $(B_n)_{n \in \N}$ is a sequence of Borel complete sections in $\Free(2^\G)$ and $(F_n)_{n \in \N}$ is a sequence of finite subsets of $\G$ such that every finite subset $F \subset \G$ is contained in $F_n$ for some $n \in \N$. 
        Then there exists a point $x \in \Free(2^\G)$ satisfying $x \in F_n \cdot B_n$ 
        for infinitely many $n \in \N$.
    \end{theo}

    Working in a free $\set{F_n}_{n \in \N}$-minimal subshift, we are able to generalize Theorem~\ref{theo:GJKS1} by removing all constraints on the shape of the sets $F_n$---we only need $|F_n|$ to be unbounded, so that Corollary~\ref{corl:subshift} may be applied to the family $\set{F_n}_{n \in \N}$:

    \begin{tcolorbox}
    \begin{theo}\label{theo:DST1}
        Suppose $(B_n)_{n \in \N}$ is a sequence of Borel complete sections in $\Free(2^\G)$ and $(F_n)_{n \in \N}$ is a sequence of finite subsets of $\G$ such that $\sup_{n \in \N}|F_n| = \infty$. Then there exists a point $x \in \Free(2^\G)$ satisfying $x \in F_n \cdot B_n$ for infinitely many $n \in \N$.
    \end{theo}
    \end{tcolorbox}

    Here is another example. Following the terminology of \cite{LLLDyn3}, we say that a set $A \subseteq \Free(2^\G)$ \emphd{traps} a point $x \in \Free(2^\G)$ (or that $x$ is \emphd{trapped} in $A$) if $\G \cdot x \subseteq A$. Gao \emph{et al.}~showed that for every Borel complete section $B \subseteq \Free(2^\G)$, some point is trapped in the union of finitely many shifts of $B$:

    \begin{theo}[{Gao--Jackson--Krohne--Seward \cite[Thm.~1.2]{GJKS_Borel}}]\label{theo:GJKS}
        If $B \subseteq \Free(2^\G)$ is a Borel complete section, then there exists a finite set $F \subset \G$ such that the set $F \cdot B$ traps a point.
    \end{theo}

    An interesting question is what additional properties the set $F$ in Theorem~\ref{theo:GJKS} may have. 
    For example, Gao \emph{et al.}~showed that for the group $\G = \Z^n$ with $n \in \N^+$, it is possible to pick $F$ so that $\|\gamma\|_1$ is odd for all $\gamma \in F$ \cite[Thm.~1.3]{GJKS_Borel}. This result is obtained with an ad hoc construction specifically tailored to the group $\Z^n$ and the ``odd $1$-norm'' property, 
    and it is unclear how far it can be generalized. 
    Corollary~\ref{corl:subshift_prime} 
    allows us to eliminate the need for such ad hoc arguments and
    completely settle the problem. 
    Namely, we show that \emph{any} large enough finite set $F \subset \G$ works: 

    \begin{tcolorbox}
    \begin{theo}\label{theo:DST}
        If $B \subseteq \Free(2^\G)$ is a Borel complete section, then there exists $n \in \N$ such that for every finite set $F \subset \G$ of size at least $n$, the set $F \cdot B$ traps a point.
    \end{theo}
    \end{tcolorbox}

    We argue by contradiction: assuming Theorem~\ref{theo:DST} fails for a Borel complete section $B \subseteq \Free(2^\G)$, we may pick sets $F_n \subset \G$ with $|F_n| = n$ so that $F_n \cdot B$ does not trap a point. We then finish the proof by working in a free $\set{F_n^{-1}}_{n \in \N}$-minimal subshift. The details are given in \S\ref{subsec:subshift_appl}.

    A version of Theorem~\ref{theo:DST} is proved in \cite[Thm.~1.2]{LLLDyn3} for \emph{open} sets $B \subseteq \Free(2^\G)$ (the proof in the open case is much simpler). We find it quite surprising that such a strong property is in fact shared by all \emph{Borel} complete sections. (Actually, our proofs of Theorems~\ref{theo:DST1} and \ref{theo:DST} do not even need the complete sections to be Borel---it is enough to assume that their intersections with every free subshift $X \subset 2^\G$ are Baire-measurable in $X$.)




    \subsection{What goes into the proof of Theorem~\ref{theo:main_prime}}\label{subsec:tools}

    Let us now say a few words about the proof of Theorem~\ref{theo:main_prime}. Although it is a result that concerns some of the most basic concepts in topological dynamics, its proof relies on a range of distinctly modern techniques. Moreover, these techniques come not only from topological dynamics, but also from descriptive set theory and combinatorics. 
    Here we describe the main streams of ideas that our argument draws upon and highlight some of the key novelties of our approach. 
    This high-level discussion will be followed in \S\ref{sec:outline} by a detailed outline of the proof of Theorem~\ref{theo:main_prime} that explains how these ideas fit together.


    \subsubsection*{Existence via genericity} Instead of constructing an explicit example of a $\G$-flow fulfilling the requirements of Theorem~\ref{theo:main_prime}, we will show that the theorem is witnessed by a ``typical'' $\G$-flow from a certain class. More precisely, we will prove that $\G$-flows with the desired properties form a dense $G_\delta$ subset in that class, and hence they are \emph{generic} in the Baire category sense.
    
    The study of generic dynamical systems is an old subject, dating back at least to the seminar work of Oxtoby and Ulam \cite{OxtobyUlam}, Halmos \cite{Halmos1,Halmos2}, and Rokhlin \cite{Rokhlin}. See \cite{Gen1,Gen2,Gen3,Gen4,Gen5,Gen6,Gen7,Gen8,Gen9} for a selection of works on genericity properties in topological dynamics from the past 25 years. 
        
    
    We borrow the idea of using genericity to obtain minimal $\G$-flows with favorable properties from the papers \cite{FTVF,FSZ}. In \cite{FTVF}, Tamuz, Vahidi Ferdowsi, and the second named author gave a genericity argument showing that every group $\G$ with a nontrivial ICC quotient admits a nontrivial proximal minimal $\G$-flow. The proof of Theorem~\ref{theo:FSZ} given in \cite{FSZ} similarly relies on genericity. The approach employed in \cite{FTVF,FSZ}, which originates in the paper \cite{FT} by Tamuz and the second named author, focuses on a specific class of $\G$-flows, namely the (closure of the) space of {strongly irreducible subshifts}. By contrast, we consider a different class that is defined specifically for the purpose of proving Theorem~\ref{theo:main_prime} (though we expect it to have other uses as well). To identify this class, we rely on the existence of a certain \emph{non-compact} action of $\G$; this is one of the central innovations of our paper. This part of the argument is described in 
    \S\ref{subsubsec:ample}.
    
    %

    \subsubsection*{Descriptive combinatorics and large-scale geometry} Our methods are strongly informed by the area called \emph{descriptive combinatorics}. This is a subject that fuses combinatorics and descriptive set theory by studying definability properties of combinatorial constructions on Polish spaces and viewing descriptive set-theoretic problems through a combinatorial lens. The systematic study of interactions between descriptive set theory and combinatorics was launched by Kechris, Solecki, and Todorcevic in their seminal paper \cite{KST} and has since developed into a rich subject with many connections to other fields, including dynamical systems. For an overview of this area, see the survey \cite{KechrisMarks} by Kechris and Marks, the introductory article \cite{Pikh_survey} by Pikhurko, and 
    the note \cite{Notices} by the first named author. 

    One of the most exciting trends in descriptive combinatorics of the past few years has been the use of concepts inspired by \emph{large-scale geometry}, pioneered by Conley, Jackson, Marks, Seward, and Tucker-Drob in their landmark paper \cite{CJMSTD}. A fundamental notion in large-scale geometry is the \emph{asymptotic dimension} of a metric space, introduced by Gromov \cite[\S1.E]{gromov1993asymptotic}. In \cite{CJMSTD}, Conley \emph{et al.}~developed a {Borel} version of asymptotic dimension and presented an array of impressive applications. They also defined a closely related parameter, the \emph{asymptotic separation index}, which turned out to be extremely useful in descriptive combinatorics \cite{ASIalgorithms,BW,BWKonig,CJMSTD,WeilacherFinDim,Grids2}. As we explain in \S\S\ref{subsubsec:example} and \ref{subsubsec:asi}, a key ingredient in our proof of Theorem~\ref{theo:main_prime} is a {continuous} variant of asymptotic separation index.
    
    Although asymptotic separation index has been associated with group actions since its initial appearance in \cite{CJMSTD}, the present work is the first to apply this notion to solve a problem of a purely dynamical nature, and we hope it will become a part of the topological dynamics toolkit.

    \subsubsection*{Definable versions of the Lov\'asz Local Lemma} The \emph{Lov\'asz Local Lemma}, or the \emph{LLL} for short, is a powerful probabilistic tool due to Erd\H{o}s and Lov\'asz \cite{EL} which has a plethora of applications throughout combinatorics; see \cite{AS,MolloyReed} for many examples. Furthermore, the LLL has recently been employed to address a number of problems in ergodic theory and topological dynamics; see, e.g., \cite{ABT,Elek,Ber_cont,LLLDyn1,LLLDyn2,LLLDyn3,LLLDyn4}. 
    
    The LLL is an existence result, and it is particularly well-suited for showing that a given structure $X$ admits a mapping $f \colon X \to k$ for some $k \in \N^+$ satisfying a specified set of constraints. Roughly speaking, in order for the LLL to apply in this context, two requirements must be met: First, a random mapping $f \colon X \to k$ should be ``likely'' to fulfill each individual constraint; second, the constraints must not interact with each other ``too much.'' (See \S\ref{subsec:LLL} for the precise statement.)


    A major research direction concerns 
    versions of the LLL that are ``constructive'' in various senses: algorithmic \cite{Beck,BGRDeterministicLLL,BMUDeterministicLLL,FG,MoserTardos}, computable \cite{RSh}, Borel/measurable \cite{CGMPT,BerDist,BerMeas,BerYuBorelLLL}, and continuous \cite{BerDist,Ber_cont}. 
    Most relevantly for the present work, Weilacher and the first named author proved a {Borel} version of the LLL under a finite asymptotic separation index assumption \cite{BW}. We employ a slight variation of their argument to establish a version of the LLL that yields \emph{continuous} maps $f \colon X \to k$ satisfying the constraints. 
    We discuss 
    the way our continuous version of the LLL becomes used in the proof of Theorem~\ref{theo:main_prime} in \S\ref{subsubsec:combi+LLL}.



    \subsubsection*{Main innovations} 
    %
    As mentioned above, the LLL has already been applied in topological dynamics in the past. However, in most previous applications, such as the ones in \cite{Elek,ABT,LLLDyn1,LLLDyn3}, it is the ``classical,'' pure existence version of the LLL that is used. By contrast, our proof of Theorem~\ref{theo:main_prime} crucially relies on a continuous variant of the LLL, thus making significant use of recent developments in descriptive combinatorics.
    
        While our ultimate goal is to find a certain continuous action of $\G$ on a \emph{compact} space, with the compactness requirement being the main challenge, our argument proceeds by first constructing an action $\G \acts X$ on a \emph{non-compact} Polish space $X$ with a technical property of being \emph{\sr}. 
        The desired $\G$-flow is then obtained by considering spaces of the form $\overline{\rho(X)}$, where $\rho \colon X \to Y$ is a continuous $\G$-equivariant map from $X$ to a $\G$-flow $Y$. 
        This somewhat counterintuitive approach separates the topological aspect of the problem---namely, compactness---from its combinatorial core, thus reducing the task to finding a Polish $\G$-space $X$ with sufficiently ``rich'' combinatorics.

        To find such a space $X$, we use asymptotic separation index. Specifically, we consider the space $\Sep(s)$ of all ``witnesses'' to the asymptotic separation index being at most $s$ and argue---using the LLL---that its free part is amply syndetic. Our work suggests that 
        the space $\Sep(s)$, which is implicit in the study of asymptotic separation index, 
        is worthy of being explored in its own right, and that its topological and dynamical properties 
         may find further use elsewhere.


    \section{Outline of the proof of Theorem~\ref{theo:main_prime}}\label{sec:outline}

    \subsection{\Sr $\G$-spaces}\label{subsubsec:ample}


    A \emphd{$\G$-space} is a topological space $X$ equipped with a continuous action $\G \acts X$. As mentioned in \S\ref{subsec:tools}, the heart of our argument is a construction of a (non-compact) Polish $\G$-space $X$ 
    satisfying the following slightly technical condition: 

    \begin{defn}[\Sr $\G$-spaces]\label{defn:amply_syndetic}
        An \emphd{\sr $\G$-space} is a 
        Polish $\G$-space $X$ 
        with the following property: 
        for every finite tuple $U_1$, \ldots, $U_k$ of nonempty open subsets of $X$, there is an integer $n \in \N$ such that for every finite set $F \subset \G$ of size at least $n$, there exists a continuous $\G$-equivariant map $\pi \colon X \to X$ for which the sets $\pi^{-1}(U_1)$, \ldots, $\pi^{-1}(U_k)$ are $F$-syndetic. 
    \end{defn}

    \begin{tcolorbox}
    \begin{theo}\label{theo:rich}
        There exists a nonempty, free, zero-dimensional 
        \sr $\G$-space.
    \end{theo}
    \end{tcolorbox}


    Recall that a topological space is \emphd{zero-dimensional} if it has a basis consisting of clopen sets. For example, the space $k^\G$ for $k \in \N^+$ and all its subspaces are zero-dimensional.
    
    Before describing the proof of Theorem~\ref{theo:rich} (which forms the bulk of the paper), let us explain how it is used to derive Theorem~\ref{theo:main_prime}. 
    %
    %
    %
    Take an arbitrary Polish $\G$-flow $Y$. (To be clear, since $Y$ is a $\G$-flow, it is a \emph{compact} Polish space.) Let $\Sub(Y)$ be the set of all subflows of $Y$. 
    We equip $\Sub(Y)$ with the \emphd{Vietoris topology}. This topology has a basis consisting of all sets of the form \label{page:Vietoris}
    \[
        \llbracket U_0; U_1, \ldots, U_k\rrbracket \,\defeq\, \set{Z \in \Sub(Y) \,:\, Z \subseteq U_0, \, Z \cap U_1 \neq \0, \, \ldots, \, Z \cap U_k \neq \0},
    \]
    where $k \in \N$ and $U_0$, \ldots, $U_k$ are open subsets of $Y$. It is a standard fact that the Vietoris topology makes $\Sub(Y)$ a compact Polish space \cite[\S4.F]{KechrisDST}. 
    Next we fix an arbitrary \sr $\G$-space $X$ and consider the following closed subset of $\Sub(Y)$:
    \[
        \Sub_X(Y) \,\defeq\, \overline{\big\{\overline{\rho(X)} \,:\, \text{$\rho \colon X \to Y$ is a continuous $\G$-equivariant map}\big\}}.
    \]
    Here horizontal lines indicate topological closure (either in $Y$ or in $\Sub(Y)$). 
    We claim that for any unbounded family $\mathcal{F}$ of finite subsets of $\G$, 
    a generic member of $\Sub_X(Y)$ is $\mathcal{F}$-minimal:

    \begin{tcolorbox}
    \begin{theo}\label{theo:generic}
        Let $X$ be an \sr $\G$-space and let $Y$ be a Polish $\G$-flow. For each 
        unbounded family $\mathcal{F}$ of finite subsets of $\G$, the set \[\set{Z \in \Sub_X(Y) \,:\, \text{$Z$ is $\mathcal{F}$-minimal}}\] is dense and $G_\delta$ in $\Sub_X(Y)$.
    \end{theo}
    \end{tcolorbox}

    With this, Theorem~\ref{theo:main_prime} easily follows:

    \begin{scproof}[ of Theorem~\ref{theo:main_prime}]
        Fix a sequence $(\mathcal{F}_n)_{n \in \N}$ of unbounded families of finite subsets of $\G$ 
        and let $X$ be a nonempty, free, zero-dimensional \sr $\G$-space given by Theorem~\ref{theo:rich}. We can find a free Polish $\G$-flow $Y$ that admits a continuous $\G$-equivariant map $\rho \colon X \to Y$, and hence $\Sub_X(Y) \neq \0$ (this is the only place where we use that $X$ is zero-dimensional; see Proposition~\ref{prop:map_to_free_flow}). Every subflow of $Y$ is free, so it suffices to argue that there exists a subflow $Z \in \Sub_X(Y)$ that is $\mathcal{F}_n$-minimal for all $n \in \N$. To this end, we note that, by Theorem~\ref{theo:generic}, \[\set{Z \in \Sub_X(Y) \,:\, \text{$Z$ is $\mathcal{F}_n$-minimal}}\] is a dense $G_\delta$ subset of $\Sub_X(Y)$ for each $n \in \N$. By the Baire category theorem, the intersection of all these sets is also dense $G_\delta$, so, in particular, it is nonempty, as desired.
    \end{scproof}

    Theorem~\ref{theo:generic} is proved in \S\ref{sec:generic}. We now turn our attention to the proof of Theorem~\ref{theo:rich}, i.e., the existence of a nonempty, free, zero-dimensional \sr $\G$-space.

    \subsection{Separators and examples of \sr $\G$-spaces}\label{subsubsec:example}

    In this subsection we describe a family of simple concrete examples of nonempty, free, zero\=/dimensional \sr $\G$-spaces. 

    We begin with some graph-theoretic terminology. Let $G$ be a (simple, undirected) graph. For a set $U \subseteq V(G)$ of vertices, we let $G[U]$ be the subgraph of $G$ \emphd{induced} by $U$, i.e., the graph with vertex set $U$ whose adjacency relation is inherited from $G$. We say that a set $U \subseteq V(G)$ is \emphd{$G$-finite} if every connected component of the graph $G[U]$ is finite.

     Now let $\G \acts X$ be an action of $\G$. Given a set $\Phi \subseteq \G$, the \emphd{Schreier graph} $\Sch(X, \Phi)$ of this action is the graph with vertex set $X$ and edge set
    \[
        E(\Sch(X,\Phi)) \,\defeq\, \big\{\set{x, \sigma \cdot x} \,:\, x \in X \text{ and } \sigma \in \Phi \text{ such that } \sigma \cdot x \neq x\big\}.
    \]
    In the case when $X = \G$ equipped with the left multiplication action $\G \acts \G$, we write \[\Cay(\G, \Phi) \,\defeq\, \Sch(\G, \Phi)\] and call $\Cay(\G,\Phi)$ the \emphd{Cayley graph} of $\G$ corresponding to $\Phi$. Note that if an action $\G \acts X$ is free, then the Schreier graph $\Sch(X, \Phi)$ is obtained by placing a copy of $\Cay(\G,\Phi)$ onto each orbit of $\G$ in $X$. For brevity, we say that a subset $U \subseteq X$ is \emphd{$\Phi$-finite} if it is $\Sch(X,\Phi)$-finite.

    \begin{defn}[Separators]\label{defn:sep}
        Let $s \in \N^+$ and let $\Phi \subset \G$ be a finite subset. An \emphd{$(s,\Phi)$-separator} is a mapping $x \colon \G \to (s+1)$ such that 
        for all $0 \leq i \leq s$, the set $x^{-1}(s) \subseteq \G$ is
        $\Phi$-finite \ep{with respect to the left multiplication action $\G \acts \G$}. The set of all $(s,\Phi)$-separators is denoted by $\Sep(s,\Phi)$.
    \end{defn}

    \begin{figure}[t]
			\centering
			\begin{tikzpicture}
                \draw[fill=blue!15,draw=white] (-2.9,-2.9) rectangle (3.9,2.9);
                \draw[rounded corners,fill=red!15,draw=lightgray,dashed] (-2-0.25,-2.5-0.25) -- (-2-0.25,2.5+0.25) -- (3+0.25,2.5+0.25) -- (3+0.25,-2.5-0.25) -- cycle;
                \draw[rounded corners,fill=blue!15,draw=lightgray,dashed] (-1.5-0.25, -1-0.25) -- (-1.5-0.25,1.5+0.25) -- (1+0.25,1.5+0.25) -- (1+0.25,1+0.25) -- (2.5+0.25,1+0.25) -- (2.5+0.25,-1.5-0.25) -- (1.5-0.25, -1.5-0.25) -- (0-0.25,-1.5-0.25) -- (0-0.25, -1-0.25) -- cycle;
                \draw[rounded corners,fill=red!15,draw=lightgray,dashed] (-1-0.25,0-0.25) -- (0+0.25,0-0.25) -- (0+0.25,1+0.25) -- (0-0.25, 1+0.25) -- (0-0.25,0.5+0.25) -- (-0.5-0.25,0.5+0.25) -- (-0.5-0.25,0+0.25) -- (-1-0.25,0+0.25) -- cycle;
                \draw[rounded corners,fill=red!15,draw=lightgray,dashed] (1-0.25,0.5+0.25) -- (2+0.25,0.5+0.25) -- (2+0.25,-0.5-0.25) -- (1-0.25,-0.5-0.25) -- cycle;
                \draw[rounded corners,fill=blue!15,draw=lightgray,dashed] (1.5-0.25,2-0.25) -- (1.5-0.25,2+0.25) -- (2.5+0.25,2+0.25) -- (2.5+0.25,2-0.25) -- cycle;
                \draw[rounded corners,fill=blue!15,draw=lightgray,dashed] (1.5-0.25,0-0.25) -- (1.5-0.25,0+0.25) -- (1.5+0.25,0+0.25) -- (1.5+0.25,0-0.25) -- cycle;
                \draw[rounded corners,fill=blue!15,draw=lightgray,dashed] (-1-0.25,-2-0.25) -- (-1-0.25,-2+0.25) -- (-1+0.25,-2+0.25) -- (-1+0.25,-2-0.25) -- cycle;
   
                \foreach \i in {-5,...,7}
                \draw[gray] (\i/2,-2.9) -- (\i/2,2.9);

                \foreach \j in {-5,...,5}
                \draw[gray] (-2.9,\j/2) -- (3.9,\j/2);


                \filldraw[blue] (-2.5,2.5) circle (2pt);
                \filldraw[red] (-2,2.5) circle (2pt);
                \filldraw[red] (-1.5,2.5) circle (2pt);
                \filldraw[red] (-1,2.5) circle (2pt);
                \filldraw[red] (-.5,2.5) circle (2pt);
                \filldraw[red] (0,2.5) circle (2pt);
                \filldraw[red] (.5,2.5) circle (2pt);
                \filldraw[red] (1,2.5) circle (2pt);
                \filldraw[red] (1.5,2.5) circle (2pt);
                \filldraw[red] (2,2.5) circle (2pt);
                \filldraw[red] (2.5,2.5) circle (2pt);
                \filldraw[red] (3,2.5) circle (2pt);
                \filldraw[blue] (3.5,2.5) circle (2pt);

                \filldraw[blue] (-2.5,2) circle (2pt);
                \filldraw[red] (-2,2) circle (2pt);
                \filldraw[red] (-1.5,2) circle (2pt);
                \filldraw[red] (-1,2) circle (2pt);
                \filldraw[red] (-.5,2) circle (2pt);
                \filldraw[red] (0,2) circle (2pt);
                \filldraw[red] (.5,2) circle (2pt);
                \filldraw[red] (1,2) circle (2pt);
                \filldraw[blue] (1.5,2) circle (2pt);
                \filldraw[blue] (2,2) circle (2pt);
                \filldraw[blue] (2.5,2) circle (2pt);
                \filldraw[red] (3,2) circle (2pt);
                \filldraw[blue] (3.5,2) circle (2pt);

                \filldraw[blue] (-2.5,1.5) circle (2pt);
                \filldraw[red] (-2,1.5) circle (2pt);
                \filldraw[blue] (-1.5,1.5) circle (2pt);
                \filldraw[blue] (-1,1.5) circle (2pt);
                \filldraw[blue] (-.5,1.5) circle (2pt);
                \filldraw[blue] (0,1.5) circle (2pt);
                \filldraw[blue] (.5,1.5) circle (2pt);
                \filldraw[blue] (1,1.5) circle (2pt);
                \filldraw[red] (1.5,1.5) circle (2pt);
                \filldraw[red] (2,1.5) circle (2pt);
                \filldraw[red] (2.5,1.5) circle (2pt);
                \filldraw[red] (3,1.5) circle (2pt);
                \filldraw[blue] (3.5,1.5) circle (2pt);

                \filldraw[blue] (-2.5,1) circle (2pt);
                \filldraw[red] (-2,1) circle (2pt);
                \filldraw[blue] (-1.5,1) circle (2pt);
                \filldraw[blue] (-1,1) circle (2pt);
                \filldraw[blue] (-.5,1) circle (2pt);
                \filldraw[red] (0,1) circle (2pt);
                \filldraw[blue] (.5,1) circle (2pt);
                \filldraw[blue] (1,1) circle (2pt);
                \filldraw[blue] (1.5,1) circle (2pt);
                \filldraw[blue] (2,1) circle (2pt);
                \filldraw[blue] (2.5,1) circle (2pt);
                \filldraw[red] (3,1) circle (2pt);
                \filldraw[blue] (3.5,1) circle (2pt);

                \filldraw[blue] (-2.5,.5) circle (2pt);
                \filldraw[red] (-2,.5) circle (2pt);
                \filldraw[blue] (-1.5,.5) circle (2pt);
                \filldraw[blue] (-1,.5) circle (2pt);
                \filldraw[red] (-.5,.5) circle (2pt);
                \filldraw[red] (0,.5) circle (2pt);
                \filldraw[blue] (.5,.5) circle (2pt);
                \filldraw[red] (1,.5) circle (2pt);
                \filldraw[red] (1.5,.5) circle (2pt);
                \filldraw[red] (2,.5) circle (2pt);
                \filldraw[blue] (2.5,.5) circle (2pt);
                \filldraw[red] (3,.5) circle (2pt);
                \filldraw[blue] (3.5,.5) circle (2pt);

                \filldraw[blue] (-2.5,0) circle (2pt);
                \filldraw[red] (-2,0) circle (2pt);
                \filldraw[blue] (-1.5,0) circle (2pt);
                \filldraw[red] (-1,0) circle (2pt);
                \filldraw[red] (-.5,0) circle (2pt);
                \filldraw[red] (0,0) circle (2pt);
                \filldraw[blue] (.5,0) circle (2pt);
                \filldraw[red] (1,0) circle (2pt);
                \filldraw[blue] (1.5,0) circle (2pt);
                \filldraw[red] (2,0) circle (2pt);
                \filldraw[blue] (2.5,0) circle (2pt);
                \filldraw[red] (3,0) circle (2pt);
                \filldraw[blue] (3.5,0) circle (2pt);

                \filldraw[blue] (-2.5,-.5) circle (2pt);
                \filldraw[red] (-2,-.5) circle (2pt);
                \filldraw[blue] (-1.5,-.5) circle (2pt);
                \filldraw[blue] (-1,-.5) circle (2pt);
                \filldraw[blue] (-.5,-.5) circle (2pt);
                \filldraw[blue] (0,-.5) circle (2pt);
                \filldraw[blue] (.5,-.5) circle (2pt);
                \filldraw[red] (1,-.5) circle (2pt);
                \filldraw[red] (1.5,-.5) circle (2pt);
                \filldraw[red] (2,-.5) circle (2pt);
                \filldraw[blue] (2.5,-.5) circle (2pt);
                \filldraw[red] (3,-.5) circle (2pt);
                \filldraw[blue] (3.5,-.5) circle (2pt);

                \filldraw[blue] (-2.5,-1) circle (2pt);
                \filldraw[red] (-2,-1) circle (2pt);
                \filldraw[blue] (-1.5,-1) circle (2pt);
                \filldraw[blue] (-1,-1) circle (2pt);
                \filldraw[blue] (-.5,-1) circle (2pt);
                \filldraw[blue] (0,-1) circle (2pt);
                \filldraw[blue] (.5,-1) circle (2pt);
                \filldraw[blue] (1,-1) circle (2pt);
                \filldraw[blue] (1.5,-1) circle (2pt);
                \filldraw[blue] (2,-1) circle (2pt);
                \filldraw[blue] (2.5,-1) circle (2pt);
                \filldraw[red] (3,-1) circle (2pt);
                \filldraw[blue] (3.5,-1) circle (2pt);

                \filldraw[blue] (-2.5,-1.5) circle (2pt);
                \filldraw[red] (-2,-1.5) circle (2pt);
                \filldraw[red] (-1.5,-1.5) circle (2pt);
                \filldraw[red] (-1,-1.5) circle (2pt);
                \filldraw[red] (-.5,-1.5) circle (2pt);
                \filldraw[blue] (0,-1.5) circle (2pt);
                \filldraw[blue] (.5,-1.5) circle (2pt);
                \filldraw[blue] (1,-1.5) circle (2pt);
                \filldraw[blue] (1.5,-1.5) circle (2pt);
                \filldraw[blue] (2,-1.5) circle (2pt);
                \filldraw[blue] (2.5,-1.5) circle (2pt);
                \filldraw[red] (3,-1.5) circle (2pt);
                \filldraw[blue] (3.5,-1.5) circle (2pt);

                \filldraw[blue] (-2.5,-2) circle (2pt);
                \filldraw[red] (-2,-2) circle (2pt);
                \filldraw[red] (-1.5,-2) circle (2pt);
                \filldraw[blue] (-1,-2) circle (2pt);
                \filldraw[red] (-.5,-2) circle (2pt);
                \filldraw[red] (0,-2) circle (2pt);
                \filldraw[red] (.5,-2) circle (2pt);
                \filldraw[red] (1,-2) circle (2pt);
                \filldraw[red] (1.5,-2) circle (2pt);
                \filldraw[red] (2,-2) circle (2pt);
                \filldraw[red] (2.5,-2) circle (2pt);
                \filldraw[red] (3,-2) circle (2pt);
                \filldraw[blue] (3.5,-2) circle (2pt);

                \filldraw[blue] (-2.5,-2.5) circle (2pt);
                \filldraw[red] (-2,-2.5) circle (2pt);
                \filldraw[red] (-1.5,-2.5) circle (2pt);
                \filldraw[red] (-1,-2.5) circle (2pt);
                \filldraw[red] (-.5,-2.5) circle (2pt);
                \filldraw[red] (0,-2.5) circle (2pt);
                \filldraw[red] (.5,-2.5) circle (2pt);
                \filldraw[red] (1,-2.5) circle (2pt);
                \filldraw[red] (1.5,-2.5) circle (2pt);
                \filldraw[red] (2,-2.5) circle (2pt);
                \filldraw[red] (2.5,-2.5) circle (2pt);
                \filldraw[red] (3,-2.5) circle (2pt);
                \filldraw[blue] (3.5,-2.5) circle (2pt);
			\end{tikzpicture}
        \caption{A piece of a $(1, \Phi)$-separator $x \colon \Z^2 \to 2$ in the group $\Z^2$, where $\Phi \defeq \set{(1,0),(0,1)}$ is the standard generating set. Here each element $\gamma \in \Z^2$ is labeled red or blue depending on the value $x(\gamma)$, splitting the Cayley graph $\Cay(\Z^2, \Phi)$ into finite monochromatic components.  
        }\label{fig:separator}
	\end{figure}
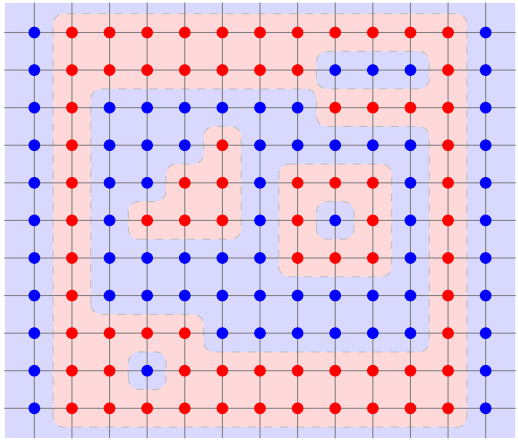

    Definition~\ref{defn:sep} is illustrated in Figure~\ref{fig:separator}. It is easy to see that for $s \in \N^+$ and finite $\Phi \subset \G$, $\Sep(s,\Phi)$ is a shift-invariant dense $G_\delta$ subset of $(s+1)^\G$. In particular, $\Sep(s,\Phi)$ is a nonempty zero-dimensional Polish space (in the subspace topology) equipped with the shift action $\G \acts \Sep(s,\Phi)$. 

    \begin{defn}[Asymptotic separators]
    For the remainder of the paper, we fix an enumeration $(\Phi_n)_{n \in \N}$ of the finite subsets of $\G$ in which every set appears infinitely often. Given $s \in \N^+$, we let
    \[
        \Sep(s) \,\defeq\, \prod_{n \in \N} \Sep(s, \Phi_n).
    \]
    The space $\Sep(s)$ is equipped with the product topology and the diagonal action of $\G$. We call the elements of $\Sep(s)$ \emphd{asymptotic $s$-separators}.
    \end{defn}
    
    Note that the action $\G \acts \Sep(s)$ is often free:

    \begin{prop}\label{prop:torsion-free}
        If $\G$ is torsion-free, then for all $s \in \N^+$, the action $\G \acts \Sep(s)$ is free.
    \end{prop}
    \begin{scproof}
        Take any $\gamma \in \G \setminus \set{\iden}$ and suppose that $\gamma \cdot x = x$ for some $x = (x_n)_{n \in \N} \in \Sep(s)$. Pick $n \in \N$ such that $\gamma \in \Phi_n$ and let $i \defeq x_n(\iden)$. Since $\gamma \cdot x_n = x_n$, we have $x_n(\gamma^m) = i$ for all $m \in \Z$. But then the infinite set $\langle \gamma \rangle$ is contained in a single component of the subgraph of $\Cay(\G, \Phi_n)$ induced by $x_n^{-1}(i)$, which contradicts the fact that $x_n \in \Sep(s, \Phi_n)$. 
    \end{scproof}

    In general, the {free part} $\Free(\Sep(s))$ of $\Sep(s)$, i.e., the set of all points $x \in \Sep(s)$ with trivial stabilizers, is a dense $G_\delta$ subset of $\Sep(s)$, because both $\Sep(s)$ and the free part of $((s+1)^\G)^\N$ are dense $G_\delta$ subsets of $((s+1)^\G)^\N$. Hence, $\Free(\Sep(s))$ is a nonempty zero-dimensional Polish space, and, by definition, the group $\G$ acts on $\Free(\Sep(s))$ freely. We will show that $\Free(\Sep(s))$ is also \sr, and therefore 
    it is an example confirming Theorem~\ref{theo:rich}:

    \begin{tcolorbox}
    \begin{theo}\label{theo:example}
        For any $s \in \N^+$, $\Free(\Sep(s))$ is an \sr $\G$-space.
    \end{theo}
    \end{tcolorbox}

    It is clear from the above discussion that Theorem~\ref{theo:example} implies Theorem~\ref{theo:rich}. In order to explain how Theorem~\ref{theo:example} is proved, we need to say a few words about the relationship between the space $\Sep(s)$ of asymptotic $s$-separators and asymptotic separation index. 

    \subsection{Asymptotic separation index and its continuous version}\label{subsubsec:asi}

    For simplicity, we will introduce asymptotic separation index for free group actions, although it is not the most general context in which it can be defined; for example, the original definition in \cite{CJMSTD} is given for so-called \emph{Borel extended metric spaces}.

    \begin{defn}[{Asymptotic separation index~\cite[Defn.~3.2]{CJMSTD}}]\label{defn:asi}
        Let $\G \acts X$ be a free Borel action of $\G$ on a Polish space $X$. The \emphd{asymptotic separation index} of $X$, in symbols $\asi(X)$, is the smallest integer $s \in \N$ such that for every finite set $\Phi \subset \G$, there exists a partition $X = X_0 \sqcup \ldots \sqcup X_s$ into $\Phi$-finite Borel sets. If there is no such $s \in \N$, we let $\asi(X) \defeq \infty$.
    \end{defn}

    We need 
    a variant of Definition~\ref{defn:asi} in which the sets $X_0$, \ldots, $X_s$ are 
    {clopen} rather than Borel:

    \begin{defn}[Continuous asymptotic separation index]\label{defn:casi}
        Let $X$ be a free zero-dimensional Polish $\G$-space. 
        The \emphd{continuous asymptotic separation index} of $X$, in symbols $\casi(X)$, is the smallest integer $s \in \N$ such that for every finite set $\Phi \subset \G$, there exists a partition $X = X_0 \sqcup \ldots \sqcup X_s$ into $\Phi$-finite clopen sets. If there is no such $s \in \N$, we let $\casi(X) \defeq \infty$.
    \end{defn}

    We can equivalently describe the (continuous) asymptotic separation index of $X$ by relating $X$ to the spaces of asymptotic $s$-separators. To this end, note that for any action $\G \acts X$ and $k \in \N^+$, we have a one-to-one correspondence
    \[
	    \big\{\text{functions $X \to k$}\big\} \quad \longleftrightarrow \quad \big\{\text{$\G$-equivariant maps $X \to k^\G$}\big\}.
	\]
	Namely, each function $f \colon X \to k$ gives rise to the $\G$-equivariant \emphd{coding map} $\pi_f \colon X \to k^\G$ via
	\[
	    \big(\pi_f(x)\big)(\gamma) \,\defeq\, f(\gamma \cdot x) \quad \text{for all } x \in X \text{ and } \gamma \in \G.
	\]
	Conversely, if $\pi \colon X \to k^\G$ is $\G$-equivariant, then the function 
    $f \colon X \to k$ such that $\pi = \pi_f$ is given by \[f(x) \,\defeq\, (\pi(x))(\iden).\]  
    Using this correspondence, we obtain the following:

    \begin{prop}[Dynamical view of continuous asymptotic separation index]\label{prop:map_to_Sep}
        Let $X$ be a free zero-dimensional Polish $\G$-space. 
        The following statements are equivalent for all $s \in \N^+$:
        \begin{enumerate}[label=\ep{\normalfont\roman*}]
            \item\label{item:casi} $\casi(X) \leq s$,
            \item\label{item:term} for every finite set $\Phi \subset \G$, there exists a continuous $\G$-equivariant map $X \to \Sep(s,\Phi)$,
            \item\label{item:product} there exists a continuous $\G$-equivariant map $X \to \Sep(s)$.
        \end{enumerate}
    \end{prop}
    \begin{scproof}
        The equivalence \ref{item:term} $\Longleftrightarrow$ \ref{item:product} is clear from the definition of $\Sep(s)$. Given a finite set $\Phi \subset \G$ and a partition $X = X_0 \sqcup \ldots \sqcup X_s$ of $X$ into $\Phi$-finite clopen sets, we obtain a continuous $\G$-equivariant map $\pi \colon X \to \Sep(s,\Phi)$ by setting $f(x) \defeq i$ for all $x \in X_i$ and taking $\pi \defeq \pi_f$. This proves \ref{item:casi} $\Longrightarrow$ \ref{item:term}. Conversely, given a continuous $\G$-equivariant map $\pi \colon X \to \Sep(s,\Phi)$, we can form a partition of $X$ into $\Phi$-finite clopen sets $X_0$, \ldots, $X_s$ by letting $X_i \defeq \set{x \in X \,:\, (\pi(x))(\iden) = i}$.
    \end{scproof}

    In particular, $\casi(\Free(\Sep(s))) \leq s$ since the inclusion map $\Free(\Sep(s)) \hookrightarrow \Sep(s)$ is continuous and $\G$-equivariant. The same proof gives a version of Proposition~\ref{prop:map_to_Sep} for the ordinary asymptotic separation index, where the word ``continuous'' is replaced by ``Borel.'' 


    In applications, the important distinction is usually between actions $\G \acts X$ with finite asymptotic separation index and those for which $\asi(X) = \infty$.\footnote{As an aside, we remark that there are currently no known examples with $1 < \asi(X) < \infty$ \cite[3191]{CJMSTD}.} In particular, it turns out that when $\asi(X) < \infty$, many otherwise intractable combinatorial problems can be solved on $X$ in a Borel way \cite{ASIalgorithms,BW,BWKonig,CJMSTD,WeilacherFinDim}. Our insight is that 
    upgrading the assumption from $\asi(X) < \infty$ to $\casi(X) < \infty$ may yield a {continuous} solution in place of a Borel one. In the next subsection we shall explain how this idea is used to find continuous $\G$-equivariant maps \[\pi \colon \Free(\Sep(s)) \to \Free(\Sep(s))\] witnessing that $\Free(\Sep(s))$ is \sr.

    \subsection{Combinatorial core of the problem and the role of the Lov\'asz Local Lemma}\label{subsubsec:combi+LLL}

    Suppose $X$ is a free zero-dimensional Polish $\G$-space. 
    By Proposition~\ref{prop:map_to_Sep}, if $\casi(X) \leq s$, then there exists a continuous $\G$-equivariant map $\pi \colon X \to \Sep(s)$. The technical substance of Theorem~\ref{theo:example} is in the following lemma, which asserts that the map $\pi$ can additionally be chosen to have properties analogous to those in Definition~\ref{defn:amply_syndetic}: 
    
    \begin{lemma}\label{lemma:key}
        For every $s \in \N^+$ and a finite tuple $U_1$, \ldots, $U_k \subseteq ((s+1)^\G)^\N$ of nonempty open sets, there is $n = n(U_1,\ldots, U_k) \in \N$ such that for all finite $F \subset \G$ with $|F| \geq n$, the following holds:
        
        \smallskip
        
        \noindent Every free zero-dimensional Polish $\G$-space $X$ 
        such that $\casi(X) \leq s$ admits a continuous $\G$-equivariant map $\pi \colon X \to \Sep(s)$ such that the sets $\pi^{-1}(U_1)$, \ldots, $\pi^{-1}(U_k)$ are $F$-syndetic.
        
    \end{lemma}

    \begin{remk}
        Note that
        in the statement of Lemma~\ref{lemma:key}, $U_1$, \ldots, $U_k$ are taken to be
        nonempty open subsets of $((s+1)^\G)^\N$ rather than of $\Sep(s)$. This is done purely for convenience and makes no difference to the content of the lemma because $\Sep(s)$ is dense in $((s+1)^\G)^\N$.
    \end{remk}

    Lemma~\ref{lemma:key} takes us most of the way toward proving Theorem~\ref{theo:example}. Indeed, suppose we are given 
    nonempty open sets $U_1$, \ldots, $U_k  \subseteq  ((s+1)^\G)^\N$ and let $n= n(U_1,\ldots, U_k) \in \N$ be the integer produced by Lemma~\ref{lemma:key}. Since $\casi(\Free(\Sep(s))) \leq s$, 
    for each finite set $F \subset \G$ of size at least $n$, there exists a continuous $\G$-equivariant map $\pi \colon \Free(\Sep(s)) \to \Sep(s)$ such that the sets $\pi^{-1}(U_1)$, \ldots, $\pi^{-1}(U_k)$ are $F$-syndetic.
    When $\G$ is torsion-free, this shows that the space $\Free(\Sep(s)) = \Sep(s)$ is \sr, as desired. In general, we need to adjust the map $\pi$ to ensure that its range is included in the free part of $\Sep(s)$. This adjustment is performed in \S\ref{sec:final}.

    To prove Lemma~\ref{lemma:key}, given a large finite set $F \subset \G$, we first 
    judiciously choose $F$-syndetic
    subsets $A_1$, \ldots, $A_k \subseteq X$ 
        and then construct a continuous $\G$-equivariant map $\pi \colon X \to \Sep(s)$ such that $\pi^{-1}(U_i) \supseteq A_i$ for all $1 \leq i \leq k$. To implement this strategy, we need to answer the question:
    \begin{tcolorbox}[width=\linewidth, sharp corners=all, colback=white!95!black, boxrule=0pt,colframe=white] 
        Under what assumptions on subsets $A_1$, \ldots, $A_k \subseteq X$ can we guarantee that there exists a continuous $\G$-equivariant map $\pi \colon X \to \Sep(s)$  with $\pi(x) \in U_i$ for all $x \in A_i$?
    \end{tcolorbox}
    We show it is enough for $A_1$, \ldots, $A_k$ to be sufficiently ``spaced out,'' in the following sense:

    \begin{defn}[$\Phi$-spaced sets]
        Given an action $\G \acts X$ and a finite set $\Phi \subset \G$, we say that a subset $A \subseteq X$ is \emphd{$\Phi$-spaced} if $(\Phi \cdot x) \cap (\Phi \cdot y) = \0$ for all distinct $x$, $y \in A$. 
    \end{defn}

    \begin{lemma}\label{lemma:from_separation_to_map}
       For every $s \in \N^+$ and a finite tuple $U_1$, \ldots, $U_k \subseteq ((s+1)^\G)^\N$ of nonempty open sets, there is a finite set $\Phi = \Phi(U_1, \ldots, U_k) \subset \G$ with the following property:

       \smallskip
       
       \noindent Suppose $X$ is a free zero-dimensional Polish $\G$-space such that $\casi(X) \leq s$, and let $A_1$, \ldots, $A_k \subseteq X$ be disjoint clopen sets whose union $A \defeq A_1 \sqcup \ldots \sqcup A_k$ is $\Phi$-spaced. Then there exists a continuous $\G$-equivariant map $\pi \colon X \to \Sep(s)$ such that $\pi(x) \in U_i$
        for all $1 \leq i \leq k$ and $x \in A_i$.
    \end{lemma}

    Lemma~\ref{lemma:from_separation_to_map} is proved in \S\ref{sec:spaced}, where the desired map $\pi$ is built explicitly. Assuming Lemma~\ref{lemma:from_separation_to_map}, the following suffices to complete the proof of Lemma~\ref{lemma:key}:

    \begin{lemma}\label{lemma:sparse_syndetic}
        For any $s$, $k \in \N^+$ and a finite set $\Phi \subset \G$, there is $n = n(s,k, \Phi) \in \N$ such that for all finite $F \subset \G$ with $|F| \geq n$, the following holds:

        \smallskip
        
        \noindent For every free zero-dimensional Polish $\G$-space $X$ such that $\casi(X) \leq s$, there exist disjoint $F$-syndetic clopen subsets $A_1$, \ldots, $A_k \subseteq X$ whose union 
          $A \defeq A_1 \sqcup \ldots \sqcup A_k$ is $\Phi$-spaced.
    \end{lemma}

    \begin{scproof}[ of Lemma~\ref{lemma:key}]
        Let $s$ and $U_1$, \ldots, $U_k$ be as in Lemma~\ref{lemma:key}. Let $\Phi = \Phi(U_1,\ldots, U_k) \subset \G$ be the finite set given by Lemma~\ref{lemma:from_separation_to_map} and let $n = n(s,k, \Phi) \in \N$ be given by  Lemma~\ref{lemma:sparse_syndetic}. Take any finite set $F \subset \G$ with $|F| \geq n$ and an arbitrary free zero-dimensional Polish $\G$-space $X$ such that $\casi(X) \leq s$. By Lemma~\ref{lemma:sparse_syndetic}, there exist disjoint $F$-syndetic clopen subsets $A_1$, \ldots, $A_k \subseteq X$ whose union $A \defeq A_1 \sqcup \ldots \sqcup A_k$ is $\Phi$-spaced. 
        By Lemma~\ref{lemma:from_separation_to_map}, we have a continuous $\G$-equivariant map $\pi \colon X \to \Sep(s)$ with $\pi(x) \in U_i$ for all $x \in A_i$, i.e., 
        $\pi^{-1}(U_i) \supseteq A_i$, which is $F$-syndetic, as desired.
    \end{scproof}

    \begin{remk}
        Lemma~\ref{lemma:sparse_syndetic} may badly fail without the assumption $\casi(X) < \infty$. For example, let $\mathbb{F}_2$ be the free group with generators $a$, $b$ and let $X$ be the free part of the shift action $\mathbb{F}_2 \acts 2^{\mathbb{F}_2}$. Then, for $\Phi \defeq \set{\iden, a}$ and for any finite set $F \subset \langle b \rangle$ (no matter how large), there is no Borel---let alone clopen---set $A \subseteq X$ that is both $F$-syndetic and $\Phi$-spaced. This is a consequence of a result of Marks \cite[Thm.~1.6]{Marks}, who showed that 
        for any Borel set $A \subseteq X$, there exists an element $x \in X$ such that either $\langle a \rangle \cdot x \subseteq A$ or $(\langle b \rangle \cdot x) \cap A = \0$.
        %
        %
    \end{remk}

    The proof of Lemma~\ref{lemma:sparse_syndetic} is where the Lov\'asz Local Lemma comes into play. 
    Consider the following purely combinatorial fact (which is implied by Lemma~\ref{lemma:sparse_syndetic}):

    \begin{quote}
        \textsl{For all finite $\Phi \subset \G$ and sufficiently large finite $F \subset \G$, 
        there is a set $A \subseteq \G$ that is both $F$-syndetic and $\Phi$-spaced \ep{with respect to the left multiplication action $\G \acts \G$}.}
    \end{quote}

    \noindent This statement has a routine proof using the LLL, and as far as we are aware, there is no other known way to establish it in full generality. This means that, to prove Lemma~\ref{lemma:sparse_syndetic}, we must rely on some version of the LLL; furthermore, this version must be able to produce \emph{clopen} sets with the desired combinatorial properties. As mentioned previously, Weilacher and the first named author \cite{BW} gave a {Borel} version of the LLL that can be used in the context of an action $\G \acts X$ with $\asi(X) < \infty$. In \S\ref{sec:CLLL}, we observe that essentially the same argument yields a {continuous} version of the LLL under the condition $\casi(X) < \infty$, which can then be applied to deduce Lemma~\ref{lemma:sparse_syndetic}.

    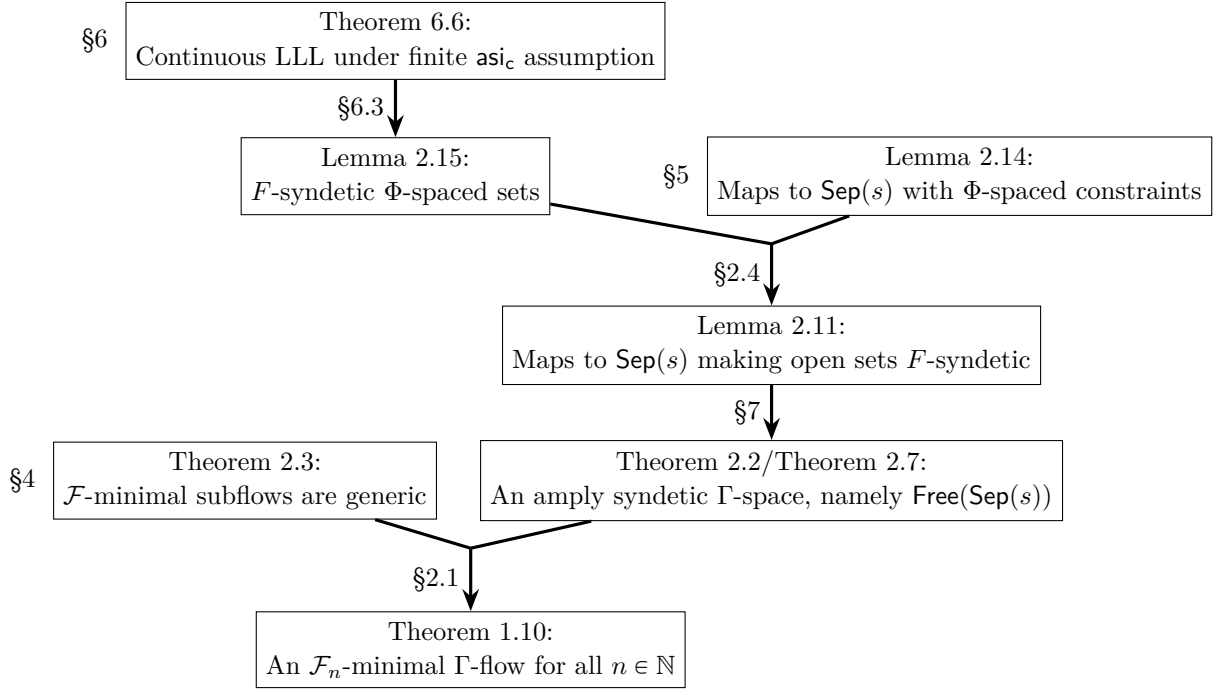
\begin{figure}[t]
			\centering
			\begin{tikzpicture}[yscale=0.9]
                    \node[draw, rectangle, align=center] (MainThm) at (0,0) {\small Theorem~\ref{theo:main_prime}: \\ \small An $\mathcal{F}_n$-minimal $\G$-flow for all $n \in \N$};



                    \node[draw, rectangle, align=center] (Generic) at (-3,2.5) {\small Theorem~\ref{theo:generic}: \\ \small $\mathcal{F}$-minimal subflows are generic};

                    \node[draw, rectangle, align=center] (Ample) at (4,2.5) {\small Theorem~\ref{theo:rich}/Theorem~\ref{theo:example}: \\ \small An \sr $\G$-space, namely $\Free(\Sep(s))$};

                    \draw [very thick] (Generic) to (0,1.5);
                    \draw [very thick] (Ample) to (0,1.5);
                    \draw [-{Stealth},very thick] (0,1.5) to node[midway,anchor=east] {\small\S\ref{subsubsec:ample}}  (MainThm);

                    \node[left=0.1cm of Generic] {\small\S\ref{sec:generic}};



                    \node[draw, rectangle, align=center] (Key) at (4,4.5) {\small Lemma~\ref{lemma:key}: \\ \small Maps to $\Sep(s)$ making open sets $F$-syndetic};

                    \draw [-{Stealth},very thick] (Key) to node[midway,anchor=east] {\small\S\ref{sec:final}}  (Ample);

                    \node[draw, rectangle, align=center] (Spaced) at (6.5,7) {\small Lemma~\ref{lemma:from_separation_to_map}: \\ \small Maps to $\Sep(s)$ with $\Phi$-spaced constraints};

                    \node[draw, rectangle, align=center] (Combi) at (-1,7) {\small Lemma~\ref{lemma:sparse_syndetic}: \\ \small $F$-syndetic $\Phi$-spaced sets};

                    \draw [very thick] (Combi) to (4,6);
                    \draw [very thick] (Spaced) to (4,6);
                    \draw [-{Stealth},very thick] (4,6) to node[midway,anchor=east] {\small\S\ref{subsubsec:combi+LLL}}  (Key);

                    \node[left=0.1cm of Spaced] {\small\S\ref{sec:spaced}};

                    \node[draw, rectangle, align=center] (LLL) at (-1,9) {\small Theorem~\ref{theo:cont_LLL}: \\ \small Continuous LLL under finite $\casi$ assumption};

                    \draw [-{Stealth},very thick] (LLL) to node[midway,anchor=east] {\small\S\ref{subsec:synd_and_spaced}}  (Combi);

                    \node[left=0.1cm of LLL] {\small\S\ref{sec:CLLL}};
			\end{tikzpicture}
        \caption{A flowchart for the proof of Theorem~\ref{theo:main_prime}.}\label{fig:flowchart}
	\end{figure}

    \subsection{Summary and a road map for the remainder of the paper}

    Figure~\ref{fig:flowchart} shows the structure of the proof of Theorem~\ref{theo:main_prime} described above. The remainder of this paper is organized as follows. After some preliminary remarks in \S\ref{sec:prelims}, we prove the remaining facts needed to deduce Theorem~\ref{theo:main_prime}, namely:
    \begin{itemize}
        \item
        Theorem~\ref{theo:generic} in \S\ref{sec:generic},
        \item
        Lemma~\ref{lemma:from_separation_to_map} in \S\ref{sec:spaced}, and
        \item
        a continuous variant of the LLL together with its consequence, Lemma~\ref{lemma:sparse_syndetic}, in \S\ref{sec:CLLL}. 
    \end{itemize}
    As explained in \S\ref{subsubsec:combi+LLL}, Lemma~\ref{lemma:from_separation_to_map} and Lemma~\ref{lemma:sparse_syndetic} combined imply Lemma~\ref{lemma:key}, which we use in \S\ref{sec:final} to complete the proof of Theorem~\ref{theo:example} and hence also of Theorem~\ref{theo:rich}. As explained in \S\ref{subsubsec:ample}, Theorems~\ref{theo:rich} and \ref{theo:generic} yield our main result, Theorem~\ref{theo:main_prime}.
    
    In the remaining sections, \S\S\ref{sec:subshift}--\ref{subsec:subshift_appl}, we establish the corollaries of our main result. Namely, we prove Corollary~\ref{corl:subshift} (the existence of subshifts with minimal subdynamics) in \S\ref{sec:subshift}, apply our results to construct disjoint $\G$-flows in \S\ref{sec:disjoint}, and study Borel complete sections in \S\ref{subsec:subshift_appl}. 
    %
    %

    \section{Remarks on zero-dimensional spaces}\label{sec:prelims}

    Here we record two simple facts about zero-dimensional spaces. The first of these facts was used to derive Theorem~\ref{theo:main_prime} from Theorems~\ref{theo:rich} and \ref{theo:generic} in \S\ref{subsubsec:ample}.

    \begin{prop}\label{prop:map_to_free_flow}
        If $X$ is a free zero-dimensional Polish $\G$-space,  
        then there exists a Polish $\G$-flow $Y$ with a continuous $\G$-equivariant map $\rho \colon X \to Y$.
    \end{prop}
    \begin{scproof}
        Let $(\gamma_n)_{n \in \N}$ be an enumeration of the non-identity elements of $\G$. By \cite[Lem.~2.3]{Ber_cont}, for each $n \in \N$, there is a continuous function $f_n \colon X \to 3$ such $f_n(x) \neq f_n(\gamma_n \cdot x)$ for all $x \in X$. (This is a continuous refinement of a result of Kechris, Solecki, and Todorcevic \cite[Prop.~4.6]{KST}.) Let
        \[
            Y_n \,\defeq\, \big\{y \in 3^\G \,:\, y(\delta) \neq y(\gamma_n \delta) \text{ for all } \delta \in \G\big\}.
        \]
        Then $Y_n \subset 3^\G$ is a subshift, 
        and we have the coding map $\pi_{f_n} \colon X \to Y_n$. Let $Y \defeq \prod_{i \in \N} Y_i$ (equipped with the diagonal action of $\G$) and define a map $\rho \colon X \to Y$ by
        \[
            \rho(x) \,\defeq\, (\pi_{f_i}(x))_{i \in \N} \,\in\, Y.
        \]
        It remains to note that the action $\G \acts Y$ is free, since for each $n \in \N$ and $y = (y_i)_{i \in \N} \in Y$, 
        \[
            (\gamma_n \cdot y_n)(\iden) \,=\, y_n(\gamma_n) \,\neq\, y_n(\iden),
        \]
        and hence $\gamma_n \cdot y \neq y$. 
    \end{scproof}

    The next observation will be used in the proof of the continuous version of the LLL in \S\ref{sec:CLLL}:

    \begin{prop}\label{prop:order}
        If $X$ is a zero-dimensional Polish space, then there exists a linear order $\leq$ on $X$ such that the set $\set{(x,y) \in X^2 \,:\, x \leq y}$ is closed.
    \end{prop}
    \begin{scproof}
        Since every zero-dimensional Polish space can be embedded into the Cantor space $2^\N$ \cite[Thm.~7.8]{KechrisDST}, it is enough to consider $X = 2^\N$. On $2^\N$, the lexicographic order works. 
    \end{scproof}

    \section{Proof of Theorem~\ref{theo:generic}: $\mathcal{F}$-minimality is generic}\label{sec:generic}

    Throughout this section, we work in a fixed Polish $\G$-flow $Y$. Recall that $\Sub(Y)$ is the space of all subflows of $Y$, equipped with the Vietoris topology (defined on p.~\pageref{page:Vietoris}).  


    
    Given $U \subseteq Y$, we define the following two subsets of $\Sub(Y)$:
    \begin{align*}
        \In(U) \,&\defeq\, \set{Z \in \Sub(Y) \,:\, Z \subseteq U},\\
        \Out(U) \,&\defeq\, \set{Z \in \Sub(Y) \,:\, U \cap Z = \0} \,=\, \In(Y \setminus U).
    \end{align*}

    \begin{lemma}\label{lemma:in-n-out}
        If $U \subseteq Y$ is an open set, then $\In(U)$ is open and $\Out(U)$ is closed in $\Sub(Y)$.
    \end{lemma}
    \begin{scproof}
        Indeed, $\In(U) = \llbracket U;\,\rrbracket$ and $\Out(U)$ is the complement of $\llbracket \Sub(Y); U\rrbracket$.
    \end{scproof}
    
    Now we fix an unbounded family $\mathcal{F}$ of finite subsets of $\G$. Given a set $U \subseteq Y$, 
    we let
    \[
        \Synd(\mathcal{F},U) \,\defeq\, \set{Z \in \Sub(Y) \,:\, \text{$U \cap Z$ is $F$-syndetic in $Z$ for some $F \in \mathcal{F}$}}.
    \]

    \begin{lemma}\label{lemma:SU}
        For each open subset $U \subseteq Y$, 
        the set $\Synd(\mathcal{F},U)$ is open in $\Sub(Y)$.
    \end{lemma}
    \begin{scproof}
        Indeed, ${\Synd}(\mathcal{F},U)$ $= \bigcup_{F \in \mathcal{F}} \In(F^{-1} \cdot U)$, which is open by Lemma~\ref{lemma:in-n-out}. 
    \end{scproof}


    \begin{prop}\label{prop:Gdelta}
        The set $\Min(\mathcal{F}) \defeq \set{Z \in \Sub(Y) \,:\, \text{$Z$ is $\mathcal{F}$-minimal}}$ is $G_\delta$ in $\Sub(Y)$.
    \end{prop}
    \begin{scproof}
        Let $(U_n)_{n \in \N}$ be a countable basis for the topology on $Y$. Then $Z \in \Sub(Y)$ is $\mathcal{F}$-minimal if and only if for all $n \in \N$, either $U_n \cap Z = \0$ or $U_n \cap Z$ is $F$-syndetic in $Z$ for some $F \in \mathcal{F}$. Hence,
        \begin{equation}\label{eq:min}
            \Min(\mathcal{F}) \,=\, \bigcap_{n \in \N} \big(\Out(U_n) \,\cup\, \Synd(\mathcal{F}, U_n)\big),
        \end{equation}
        which is $G_\delta$ by Lemmas~\ref{lemma:in-n-out} and \ref{lemma:SU}.
    \end{scproof}

    We also need the following observation:

    \begin{lemma}\label{lemma:closure}
        Let $V \subseteq Y$ and suppose that $A \subseteq Y$ is a $\G$-invariant set such that $V \cap A$ is $F$-syndetic in $A$ for some finite set $F \subset \G$. Then $\overline{V} \cap \overline{A}$ is $F$-syndetic in $\overline{A}$.
    \end{lemma}
    \begin{scproof}
        By assumption, 
        $F^{-1} \cdot V \supseteq A$. Note that for each $\sigma \in F$,
        \[\sigma^{-1} \cdot \overline{V} \,=\, \overline {\sigma^{-1} \cdot V}.\]
        Since $F$ is finite, it follows that $F^{-1} \cdot \overline{V} = \overline{F^{-1} \cdot V} \supseteq \overline{A}$, as desired.
    \end{scproof}


    At this point, we fix an \sr $\G$-space $X$ and recall that
    \[
        \Sub_X(Y) \,\defeq\, \overline{\big\{\overline{\rho(X)} \,:\, \text{$\rho \colon X \to Y$ is a continuous $\G$-equivariant map}\big\}} \,\subseteq\, \Sub(Y).
    \]
    We are interested in the following subset of $\Sub_X(Y)$:
    \[
        \Min_X(\mathcal{F}) \,\defeq\, \Min(\mathcal{F}) \cap \Sub_X(Y) \,=\, \set{Z \in \Sub_X(Y) \,:\, \text{$Z$ is $\mathcal{F}$-minimal}}.
    \]
    Our goal is to show that $\Min_X(\mathcal{F})$ is a dense $G_\delta$ subset of $\Sub_X(Y)$. Since $\Min_X(\mathcal{F})$ is a $G_\delta$ set by Proposition~\ref{prop:Gdelta}, the only thing left to check is density.

    \begin{prop}\label{prop:dense}
        The set $\Min_X(\mathcal{F})$ is dense in $\Sub_X(Y)$.
    \end{prop}
    \begin{scproof}
        Thanks to formula \eqref{eq:min} and the Baire category theorem, it suffices to argue that the set
        \[
            \big(\Out(U) \,\cup\, \Synd(\mathcal{F}, U)\big) \,\cap\, \Sub_X(Y)
        \]
        is dense in $\Sub_X(Y)$ for every open $U \subseteq Y$. 
        To this end, suppose 
        \begin{equation}\label{eq:nonempty_basic}
            \llbracket U_0; U_1, \ldots, U_k \rrbracket \,\cap\, \Sub_X(Y) \,\neq \,\0
        \end{equation}
        for some basic open set $\llbracket U_0; U_1, \ldots, U_k \rrbracket$ in $\Sub(Y)$. We need to find a subflow
        \[
            Z \,\in\, \llbracket U_0; U_1, \ldots, U_k \rrbracket \,\cap\, \big(\Out(U) \,\cup\, \Synd(\mathcal{F}, U)\big) \,\cap\, \Sub_X(Y).
        \]
        By \eqref{eq:nonempty_basic} and the definition of $\Sub_X(Y)$, there is a continuous $\G$-equivariant map $\rho \colon X \to Y$ with
        \[
            \overline{\rho(X)} \,\in\, \llbracket U_0; U_1, \ldots, U_k \rrbracket.
        \]
        If $\overline{\rho(X)} \in \Out(U)$, then taking $Z \defeq \overline{\rho(X)}$ completes the proof. Otherwise, we have
        \[
            \overline{\rho(X)} \cap U_1 \neq \0, \quad \ldots, \quad \overline{\rho(X)} \cap U_k \neq \0, \quad \overline{\rho(X)} \cap U \neq \0.
        \]
        Since the sets $U_1$, \ldots, $U_k$, $U$ are open, this is equivalent to
        \[
            {\rho(X)} \cap U_1 \neq \0, \quad \ldots, \quad {\rho(X)} \cap U_k \neq \0, \quad {\rho(X)} \cap U \neq \0.
        \]
        We may pick open sets $V_1$, \ldots, $V_k$, $V$ such that $\overline{V_1} \subseteq U_1$, \ldots, $\overline{V_k} \subseteq U_k$, $\overline{V} \subseteq U$ and
        \[
            {\rho(X)} \cap V_1 \neq \0, \quad \ldots, \quad {\rho(X)} \cap V_k \neq \0, \quad {\rho(X)} \cap V \neq \0.
        \]
        Then $\rho^{-1}(V_1)$, \ldots, $\rho^{-1}(V_k)$, $\rho^{-1}(V)$ are nonempty open subsets of $X$. Since $X$ is \sr, for any sufficiently large set $F \in \mathcal{F}$ (which exists because $\mathcal{F}$ is unbounded), we have a
        %
        continuous $\G$-equivariant map $\pi \colon X \to X$ such that the sets
        \[
            (\rho \circ \pi)^{-1}(V_1), \quad \ldots, \quad (\rho \circ \pi)^{-1}(V_k), \quad (\rho \circ \pi)^{-1}(V)
        \]
        are $F$-syndetic in $X$. Equivalently, the sets
        \[
            V_1 \cap (\rho \circ \pi)(X), \quad \ldots, \quad V_k \cap (\rho \circ \pi)(X), \quad V \cap (\rho \circ \pi)(X)
        \]
        are $F$-syndetic in $(\rho \circ \pi)(X)$. 
        Now we claim that we can take 
        \[
            Z \,\defeq\, \overline{(\rho \circ \pi)(X)}. 
        \]
        Indeed, $\rho \circ \pi \colon X \to Y$ is a continuous $\G$-equivariant map, and hence $Z \in \Sub_X(Y)$. 
        Since $U \supseteq \overline{V}$, the set $U \cap Z$ is $F$-syndetic in $Z$ by Lemma~\ref{lemma:closure}, and thus $Z \in \Synd(\mathcal{F},U)$. The same reasoning shows that the sets $U_1 \cap Z$, \ldots, $U_k \cap Z$ are $F$-syndetic in $Z$ as well; in particular, they are nonempty.  Finally, $Z \subseteq \overline{\rho(X)} \subseteq U_0$ by construction, and therefore $Z \in \llbracket U_0; U_1, \ldots, U_k\rrbracket$, as desired.
    \end{scproof}

    \begin{scproof}[ of Theorem~\ref{theo:generic}]
        Follows by combining Propositions~\ref{prop:Gdelta} and \ref{prop:dense}.
    \end{scproof}

    \section{Proof of Lemma~\ref{lemma:from_separation_to_map}: Maps to $\Sep(s)$ with $\Phi$-spaced constraints}\label{sec:spaced}

    In this section we prove Lemma~\ref{lemma:from_separation_to_map}, restated here for ease of reference:

    \begin{lemmacopy}{lemma:from_separation_to_map}
       For every $s \in \N^+$ and a finite tuple $U_1$, \ldots, $U_k \subseteq ((s+1)^\G)^\N$ of nonempty open sets, there is a finite set $\Phi = \Phi(U_1, \ldots, U_k) \subset \G$ with the following property:

       \smallskip
       
       \noindent Suppose $X$ is a free zero-dimensional Polish $\G$-space such that $\casi(X) \leq s$, and let $A_1$, \ldots, $A_k \subseteq X$ be disjoint clopen sets whose union $A \defeq A_1 \sqcup \ldots \sqcup A_k$ is $\Phi$-spaced. Then there exists a continuous $\G$-equivariant map $\pi \colon X \to \Sep(s)$ such that $\pi(x) \in U_i$
        for all $1 \leq i \leq k$ and $x \in A_i$.
    \end{lemmacopy}
    \begin{scproof}
        We may assume that there exist a finite set $D \subset \G$ containing $\iden$, a natural number $N \in \N$, and mappings $\phi_{i,n} \colon D \to (s+1)$ for $1 \leq i \leq k$, $0 \leq n < N$ such that 
        $
            U_i = \prod_{n \in \N} U_{i,n} 
        $,
        where 
        \[
            U_{i,n} \,\defeq\, \begin{cases}
                \big\{x \in (s+1)^{\G} \,:\, \rest{x}{D} = \phi_{i,n} \big\} &\text{if } n < N,\\
                (s+1)^\G &\text{if } n \geq N.
            \end{cases}
        \]
        No generality is lost because sets of this type form a basis for the topology on $((s+1)^{\G})^{\N}$.
        
        Recall that $\Sep(s) = \prod_{n \in \N} \Sep(s, \Phi_n)$, where $(\Phi_n)_{n \in \N}$ is a list of all finite subsets of $\G$ that includes each set infinitely often. We let $\Phi_n^* \defeq \Phi_n \cup \Phi_n^{-1} \cup \set{\iden}$ and define
        \[
             \Phi \,\defeq\, \bigcup_{n < N} \Phi_n^* DD^{-1}\Phi_n^*.
        \]
        We claim that this finite set $\Phi$ works.
        
        Let $X$, $A_1$, \ldots, $A_k$, $A$ be as in the statement of the lemma. We seek a continuous $\G$-equivariant map $\pi \colon X \to \Sep(s)$ such that $\pi(x) \in U_i$ for all $x \in A_i$. Since $U_i = \prod_{n \in \N} U_{i,n}$, it is enough to argue that for each $n \in \N$, there exists a continuous $\G$-equivariant map
            \[
                \pi_n \colon X \to \Sep(s, \Phi_n)
        \]
        such that $\pi_n(x) \in U_{i,n}$ for all $x \in A_i$, as then we can define the desired map $\pi \colon X \to \Sep(s)$ via
        \[  
            \pi(x) \,\defeq\, (\pi_n(x))_{n \in \N}.
        \]
        Thus, for the remainder of the proof, we fix some $n \in \N$.
        
        If $n \geq N$, then 
        we can let $\pi_n \colon X \to \Sep(s,\Phi_n)$ be an arbitrary continuous $\G$-equivariant map, which exists by Proposition~\ref{prop:map_to_Sep} as $\casi(X) \leq s$.
            
            Now suppose that $n < N$. 
            Since $\casi(X) \leq s$, we may fix a continuous map $h \colon X \to (s+1)$ such that the sets $h^{-1}(0)$, \ldots, $h^{-1}(s)$ are $\Phi$-finite. 
            %
            Define
            \begin{equation}\label{eq:modify}
                f(x) \,\defeq\, \begin{cases}
                    \phi_{i,n}(\delta) &\text{if $\delta \in D$ and $\delta^{-1}\cdot x \in A_i$},\\
                    h(x) &\text{if $x \notin D \cdot A$}.
                \end{cases}
            \end{equation}
            Note that since $\Phi \supseteq D$ and $A$ is $\Phi$-spaced, for each point $x \in D \cdot A$, there exists a unique pair $(\delta,i) \in D \times \set{1, \ldots, k}$ such that $\delta^{-1}\cdot x \in A_i$, and hence formula \eqref{eq:modify} describes a well-defined continuous function $f \colon X \to (s+1)$. Informally, to obtain the function $f$, we modify $h$ by copying the mapping $\phi_{i,n} \colon D \to (s+1)$ onto each set of the form $D \cdot x$ with $x \in A_i$.
            
            We claim that the coding map $\pi_n \defeq \pi_f \colon X \to (s+1)^\G$ corresponding to $f$ has the desired properties. The definition of $f$ ensures that for all $x \in A_i$, we have $\rest{\pi_n(x)}{D} = \phi_{i,n}$, i.e., $\pi_n(x) \in U_{i,n}$. It remains to verify that $\pi_n(x) \in \Sep(s, \Phi_n)$ for all $x \in X$. 

            To this end, fix $0 \leq i \leq s$ and let $G_i$ be the subgraph of $\Sch(X, \Phi_n)$ induced by the set $f^{-1}(i)$. Similarly, let $H_i$ be the subgraph of $\Sch(X,\Phi)$ induced by $h^{-1}(i)$. 
            We need to argue that all connected components of $G_i$ are finite. 
            Suppose, toward a contradiction, that $C$ is an infinite component of $G_i$. Since the set $A$ is $\Phi$-spaced and $\Phi \supseteq \Phi_n^* D$, there is no infinite path in $\Sch(X,\Phi_n)$ consisting solely of vertices in $D \cdot A$. Hence, the set $C \setminus (D \cdot A)$ must be infinite. This, however, is impossible, as any two vertices in $C \setminus (D \cdot A)$ belong to the same \ep{finite} component of $H_i$. Indeed, suppose 
            that $(x_0,x_1,\ldots,x_\ell)$ is a shortest path in $C$ such that $x_0$, $x_\ell \in C \setminus (D \cdot A)$ belong to distinct components of $H_i$ (note that they must belong to $H_i$ because $h$ and $f$ may only differ on $D \cdot A$). As $\Phi \supseteq \Phi_n$, we have $\ell \geq 2$. By the minimality of the path, $x_1$, \ldots, $x_{\ell - 1} \in D \cdot A$, and since $A$ is $\Phi$-spaced, we conclude that $x_1$, \ldots, $x_{\ell - 1} \in D \cdot y$ for some $y \in A$. 
            But then
            \begin{align*}
                x_{\ell} \,\in\, \Phi_n^* \cdot x_{\ell - 1} \,\subseteq\, \Phi_n^*D \cdot y \,\subseteq\, \Phi_n^*DD^{-1} \cdot x_1 \,\subseteq\,  \Phi_n^*DD^{-1}\Phi_n^* \cdot x_0 \,\subseteq\, \Phi \cdot x_0,
            \end{align*}
            which implies that $x_\ell$ and $x_0$ are either equal or adjacent in $H_i$; a contradiction.
    \end{scproof}

    \section{Continuous Lov\'asz Local Lemma via continuous asymptotic separation index}\label{sec:CLLL}

    \subsection{Classical Lov\'asz Local Lemma}\label{subsec:LLL}

    Here we state the LLL in its classical combinatorial form, before moving on to the continuous setting in \S\ref{subsec:CLLL}. Given a family $\mathcal{B}$ of ``bad'' random events in a probability space, the LLL provides a sufficient condition that guarantees that all these ``bad'' events can be avoided. To formulate it, we need the following notion:

    \begin{defn}[Dependency relations]
        Let $\mathcal{B}$ be a family of random events in a probability space $(\Omega, \P)$. A \emphd{dependency relation} on $\mathcal{B}$ is a reflexive binary relation $\sim$ 
        such that every event $B \in \mathcal{B}$ is mutually independent\footnote{Recall that random events $A$, $A'$ are \emphd{independent} if $\P[A \cap A'] = \P[A] \P[A']$, and an event $A$ is \emphd{mutually independent} from a family of events $\mathcal{A}$ if it is independent from every Boolean combination of the events in $\mathcal{A}$.} from the events $\set{B' \in \mathcal{B} \,:\, B' \not\sim B}$.
    \end{defn}

    \begin{exmp}\label{exmp:vbl}
        Let $(\Omega_i, \P_i)_{i \in I}$ be a family of probability spaces and let $(\Omega, \P) \defeq \prod_{i \in I} (\Omega_i, \P_i)$ be their product space. Suppose $\mathcal{B}$ is a family of random events in $\Omega$ such that the outcome of each event $B \in \mathcal{B}$ is determined by some nonempty subset $\dom(B) \subseteq I$ of coordinates. Then the relation $\sim$ such that $B \sim B'$ if and only if $\dom(B) \cap \dom(B' )\neq \0$ is a dependency relation on $\mathcal{B}$ \cite[41]{MolloyReed}. This is by far the most common setting in which the LLL is applied.
    \end{exmp}

    \begin{theo}[{Lov\'asz Local Lemma \cites{EL}{SpencerRamsey}[Corl.~5.1.2]{AS}}]\label{theo:LLL}
        Let $\mathcal{B}$ be a finite family of random events in a probability space $(\Omega, \P)$ and let $\sim$ be a dependency relation on $\mathcal{B}$. Suppose that there exist $p \in [0,1)$ and $d \in \N^+$ such that:
        \begin{itemize}
            \item $\P[B] \leq p$ and $|\set{B' \in \mathcal{B} \,:\, B' \sim B}| \leq d$ for all $B \in \mathcal{B}$,
            \item $\mathsf{e} p d \leq 1$, where $\mathsf{e} = 2.718\ldots$ is the base of the natural logarithm. 
        \end{itemize}
        Then $\P\left[\bigcup \mathcal{B}\right] < 1$ and, in particular, $\Omega \setminus \bigcup \mathcal{B} \neq \0$.
    \end{theo}

    In practice, the ``in particular'' part of Theorem~\ref{theo:LLL} is often valid for infinite families $\mathcal{B}$ as well. For example, we have the following:

    \begin{corl}[{Infinite LLL}]\label{corl:inf_LLL}
        Suppose $\Omega$ is a compact space endowed with a Borel probability measure $\P$. Let $\mathcal{B}$ be a family of open subsets of $\Omega$ and let $\sim$ be a dependency relation on $\mathcal{B}$. If there exist $p \in [0,1)$ and $d \in \N^+$ as in Theorem~\ref{theo:LLL},  
        then $\Omega \setminus \bigcup \mathcal{B} \neq \0$.
    \end{corl}
    \begin{scproof}
        By compactness, we just need to argue that $\Omega \setminus \bigcup \mathcal{B}' \neq \0$ for every finite subfamily $\mathcal{B}' \subseteq \mathcal{B}$, which is true by Theorem~\ref{theo:LLL}.
    \end{scproof}

    Let us now describe the framework in which the LLL shall be employed in this paper. Fix $k \in \N^+$ and let $W \subset \G$ be a finite set, called a \emphd{window}. Given a set $\mathcal{P} \subseteq k^W$ of \emphd{patterns}, we define
    \[
        \Sigma(k, W, \mathcal{P}) \,\defeq\, \big\{x \in k^\G \,:\, \rest{(\gamma \cdot x)}{W} \in \mathcal{P} \text{ for all } \gamma \in \G\big\},
    \]
    where, as usual, the vertical line indicates the restriction of a function to a subset of its domain. The set $\Sigma(k, W, \mathcal{P})$ is closed and shift-invariant, so when $\Sigma(k,W,\mathcal{P}) \neq \0$, it is a subshift. Subshifts of the form $\Sigma(k, W, \mathcal{P})$ are called \emphd{subshifts of finite type}; see \cite{subshifts1,subshifts2} for more background on their role in symbolic dynamics.
    
    A simple application of the LLL provides a condition that implies $\Sigma(k, W, \mathcal{P})$ is nonempty:

    \begin{prop}\label{prop:SFT}
       Let $k \in \N^+$, let $\0\neq W \subset \G$ be a finite set, and let $\mathcal{P} \subseteq k^W$. Suppose that
       \[
            \mathsf{e} \left(1 - \frac{|\mathcal{P}|}{k^{|W|}}\right) |W|^2 \,\leq\, 1.
       \]
       Then $\Sigma(k, W, \mathcal{P}) \neq \0$.
    \end{prop}
    \begin{scproof}
        We equip the set $k = \set{0,\ldots, k-1}$ with the uniform probability measure and work in the space $k^\G$ with the product measure $\P$. For $\gamma \in \G$, define $B_\gamma \subseteq k^\G$ by
        \[
            B_\gamma \,\defeq\, \big\{x \in k^\G \,:\, \rest{(\gamma \cdot x)}{W} \notin \mathcal{P}\big\},
        \]
        and let $\mathcal{B} \defeq \set{B_\gamma \,:\, \gamma \in \G}$. Note that each set $B_\gamma$ is clopen.  
        As in Example~\ref{exmp:vbl}, we observe that the membership of a point $x \in k^\G$ in $B_\gamma$ is determined by the restriction of $x$ to $W\gamma$, and hence we can define a dependency relation $\sim$ on $\mathcal{B}$ via 
        $B_\gamma \sim B_{\gamma'}$ $\Longleftrightarrow$ 
        $(W\gamma) \cap (W\gamma') \neq \0$. Each event $B_\gamma$ satisfies
        \[
            \P[B_\gamma] \,=\, 1 - \frac{|\mathcal{P}|}{k^{|W|}} \,\eqqcolon\, p \qquad \text{and} \qquad |\set{B_{\gamma'} \in \mathcal{B} \,:\, B_\gamma \sim B_{\gamma'}}| \,=\, |W^{-1}W\gamma| \,\leq\, |W|^2 \,\eqqcolon\, d,
        \]
        so we may apply Corollary~\ref{corl:inf_LLL} to conclude that $\Sigma(k, W, \mathcal{P}) = k^\G \setminus \bigcup \mathcal{B} \neq \0$, as desired.
    \end{scproof}

    \subsection{A continuous version}\label{subsec:CLLL}

    We can now formulate the continuous variant of the LLL that we use to prove Lemma~\ref{lemma:sparse_syndetic}:

    \begin{tcolorbox}
    \begin{theo}\label{theo:cont_LLL}
        Let $s$, $k \in \N^+$, let $\0\neq W \subset \G$ be a finite set, and let $\mathcal{P} \subseteq k^W$. Suppose that
       \begin{equation}\label{eq:cont_LLL}
            \mathsf{e}^{s+1} \left(1 - \frac{|\mathcal{P}|}{k^{|W|}}\right) |W|^{2(s+1)} \,\leq\, 1.
       \end{equation}
       Then every free zero-dimensional Polish $\G$-space $X$ 
       such that $\casi(X) \leq s$ admits a continuous $\G$-equivariant map $\pi \colon X \to \Sigma(k, W, \mathcal{P})$.
    \end{theo}
    \end{tcolorbox}

    A straightforward application of a result of Weilacher and the first named author \cite[Thm.~1.29]{BW} 
    yields a {Borel} $\G$-equivariant map $\pi \colon X \to \Sigma(k,W,\mathcal{P})$ assuming $\asi(X) \leq s$. To prove Theorem~\ref{theo:cont_LLL}, we observe that when $\casi(X) \leq s$, a careful implementation of the construction from \cite{BW} 
        makes the map $\pi$ continuous.
        To make the presentation self-contained, we describe the construction here. The main tool we use is the \emph{method of conditional probabilities}, which is a common derandomization technique in computer science 
        \cites[\S16]{AS}[\S5.6]{RandAlg}. An analogous approach was employed in \cite{Ber_cont} to establish a different continuous version of the LLL, and similar arguments have appeared in other related contexts as well, for example in \cite[Thm.~3.6]{FG} by Fischer and Ghaffari.
    
    \begin{scproof}

        
        Thanks to the coding map correspondence $f \longleftrightarrow \pi_f$, our goal is to find a continuous function $f \colon X \to k$ such that for all $x \in X$, the map $W \to k$ sending each $\sigma \in W$ to $f(\sigma \cdot x)$ belongs to $\mathcal{P}$. To this end, we introduce a family of ``bad'' random events analogous to the one in the proof of Proposition~\ref{prop:SFT}. Let $k = \set{0,\ldots, k-1}$ carry the uniform probability measure and consider the 
        space $k^X$ with the product measure $\P$. For each $x \in X$, let 
        \[
            B_x \,\defeq\, \big\{f \in k^X \,:\, \text{the map $W \to k \colon \sigma \mapsto f(\sigma \cdot x)$ is not in $\mathcal{P}$}\big\},
        \]
        and set $\mathcal{B} \defeq \set{B_x \,:\, x \in X}$. Each set $B_x$ is clopen in $k^X$. The relation $\sim$ on $\mathcal{B}$ given by $B_x \sim B_y$ $\Longleftrightarrow$ $(W \cdot x) \cap (W \cdot y) \neq \0$ is a dependency relation, and every event $B_x$ satisfies
        \[
            \P[B_x] \,=\, 1 - \frac{|\mathcal{P}|}{k^{|W|}} \,\eqqcolon\, p \qquad \text{and} \qquad |\set{B_{y} \in \mathcal{B} \,:\, B_x \sim B_{y}}| \,=\, |W^{-1}W \cdot x| \,\leq\, |W|^2 \,\eqqcolon\, d.
        \]
        At this point we may apply Corollary~\ref{corl:inf_LLL} to obtain a function $f \in k^X \setminus \bigcup \mathcal{B}$. However, this is not enough for our purposes, as we also need $f$ to be {continuous}.

        To achieve this aim, we will use the assumption that $\casi(X) \leq s$. Define $\Phi \defeq WW^{-1}$, and fix a partition $X = X_0 \sqcup \ldots \sqcup X_s$ into $\Phi$-finite clopen sets. Let $G_i$ be the subgraph of the Schreier graph $\Sch(X,\Phi)$ induced by $X_i$ and let $\mathcal{C}_i$ be the set of all connected components of $G_i$. Note that all components of $G_i$ are finite. The set $\Phi$ is chosen to guarantee that each set of the form $W \cdot x$ meets at most one component of $G_i$.
        
        %
        Given a partial map $h \colon X \pto k$ and a point $x \in X$, we write $\P[B_x \,\vert\, h]$ for the conditional probability that a random function $f \in k^X$ is in $B_x$
        given that $f$ agrees with $h$ on $(W \cdot x) \cap \dom(h)$. Note that when $\dom(h) = X$, we have $\P[B_x \,\vert\, h] = 1$ if $h \in B_x$ and $\P[B_x \,\vert\, h] = 0$ otherwise. 
        Inducting on $0 \leq i \leq s$, we shall define continuous functions $f_i \colon X_i \to k$ such that for all $x \in X$,
        \[
            \P[B_x \,\vert\, f_0 \cup \ldots \cup f_i]  \,<\, (\mathsf{e} d)^{-(s - i)}.
        \]
        Then the map $f \defeq f_0 \cup \ldots \cup f_s \colon X \to k$ satisfies $\P[B_x \,\vert\, f] < 1$ for all $x \in X$. As $\dom(f) = X$, we conclude that $\P[B_x \,\vert\, f] = 0$ for all $x \in X$, 
        i.e., $f \notin \bigcup \mathcal{B}$, as desired.

        To construct the functions $f_0$, \ldots, $f_s$, we first observe that, by \eqref{eq:cont_LLL}, 
        \[
            \P[B_x] \,=\, p \,<\, (\mathsf{e} d)^{-(s+1)} \quad \text{for all } x \in B_x.
        \]
        (The inequality is strict because $p$ is rational, while $(\mathsf{e}d)^{s+1}$ is not.) Now suppose we already have continuous functions $f_0 \colon X_0 \to k$, \ldots, $f_{i-1} \colon X_{i-1} \to k$ for some $0 \leq i \leq s$ such that for all $x \in X$,
        \begin{equation}\label{eq:ind}
            \P[B_x \,\vert\, f_0 \cup \ldots \cup f_{i-1}] \,<\, (\mathsf{e}d)^{-(s-i+1)}.
        \end{equation}
        To define $f_i$, we consider the components of the graph $G_i$ individually. For $C \in \mathcal{C}_i$, we say that a function $h \colon C \to k$ is \emphd{good} if
        \[
                \P[B_x \,\vert\, f_0 \cup \ldots \cup f_{i-1} \cup h] \,<\, (\mathsf{e}d)^{-(s-i)} \quad \text{ for all } x \in X.
            \]
        
        \begin{claim}
            For each component $C \in \mathcal{C}_i$ of $G_i$, there exists a good function $h \colon C \to k$. 
        \end{claim}
        \begin{claimproof}
            Working in the space $k^C$ with the product measure $\P_C$, we consider the random events
            \[
                B_{x,C} \,\defeq\, \big\{h \in k^C \,:\, \P[B_x \,\vert\, f_0 \cup \ldots \cup f_{i-1} \cup h] \geq (\mathsf{e}d)^{-(s-i)} \big\},
            \]
            and set $\mathcal{B}_{C} \defeq \set{B_{x,C} \,:\, x \in X}$. Since each event $B_{x,C}$ is determined by the restriction of $h$ to the set $(W \cdot x) \cap C \subseteq W \cdot x$, 
            relating $B_{x,C}$ and $B_{y,C}$ if and only if $(W \cdot x) \cap (W \cdot y)$ gives a dependency relation on $\mathcal{B}_{x,C}$. 
            Using Markov's inequality and \eqref{eq:ind}, we see that for all $x \in X$,
            \[
                \P_C[B_{x,C}] \,\leq\, \frac{\P[B_x \,\vert\, f_0 \cup \ldots \cup f_{i-1}]}{(\mathsf{e}d)^{-(s-i)}} \,<\, \frac{(\mathsf{e}d)^{-(s-i+1)}}{(\mathsf{e}d)^{-(s-i)}} \,=\, \frac{1}{\mathsf{e}d}.
            \]
            It follows that $\mathcal{B}_{C}$ satisfies the assumptions of the LLL, and hence a desired function $h \in k^C \setminus \bigcup \mathcal{B}_C$ exists by Corollary~\ref{corl:inf_LLL}.
        \end{claimproof}

        Fix a closed linear order $\leq$ on $X$ given by Proposition~\ref{prop:order}. Take any $C \in \mathcal{C}_i$ and list the elements of $C$ in the increasing order as $x_1 < \cdots < x_{|C|}$. A function $h \colon C \to k$ can then be represented by the tuple $(h(x_1), \ldots, h(x_{|C|}))$, and we let $h_C \colon C \to k$ be the good function 
        for which this tuple is lexicographically minimal. 
        Define
        $
            f_i \defeq \bigcup_{C \in \mathcal{C}_i} h_C
        $.
        We claim that $f_i \colon X_i \to k$ is as desired.
        
        First, we observe that for any $x \in X$,
        \[
            \P[B_x \,\vert\, f_0 \cup \ldots \cup f_{i-1} \cup f_i] \,<\, (\mathsf{e}d)^{-(s-i)}.
        \]
        Indeed, the value $\P[B_x \,\vert\, f_0 \cup \ldots \cup f_{i-1} \cup f_i]$ is determined by the restriction of $f_i$ to $(W \cdot x) \cap X_i$. If $(W \cdot x) \cap X_i = \0$, then the desired bound holds by \eqref{eq:ind}. 
        On the other hand, if $(W \cdot x) \cap X_i \neq \0$, then there is a unique component $C \in \mathcal{C}_i$ such that $(W \cdot x) \cap X_i = (W \cdot x) \cap C$, in which case
        \[
            \P[B_x \,\vert\, f_0 \cup \ldots \cup f_{i-1} \cup f_i] \,=\, \P[B_x \,\vert\, f_0 \cup \ldots \cup f_{i-1} \cup h_C] \,<\, (\mathsf{e}d)^{-(s-i)},
        \]
        because $h_C$ is good.

        It remains to verify that $f_i$ is continuous. To this end, consider any point $x \in X_i$ and let $R_x$ be the (finite) set of all group elements $\gamma \in \G$ such that $\gamma \cdot x$ is in the same component of $G_i$ as $x$. It is evident from the definition of $f_i$ that the value $f_i(x)$ is determined by the following information:
        \begin{itemize}
            \item the set $R_x$,
            \item the restrictions of the functions $f_0$, \ldots, $f_{i-1}$ to $\Phi R_x \cdot x$,
            \item the restriction of the order $\leq$ to the set $R_x \cdot x$.
        \end{itemize}
        Since the sets $X_0$, \ldots, $X_s$ are clopen, the functions $f_0$, \ldots, $f_{i-1}$ and the action 
        $\G \acts X$ are continuous, 
        and the order $\leq$ is closed, this shows that $f_i(x) = f_i(y)$ for all $y$ in some open neighborhood of $x$, and hence $f_i$ is continuous, as desired.
     \end{scproof}

    \subsection{Application: Proof of Lemma~\ref{lemma:sparse_syndetic}}\label{subsec:synd_and_spaced}

    We are now ready to prove Lemma~\ref{lemma:sparse_syndetic}, restated here for ease of reference:

    \begin{lemmacopy}{lemma:sparse_syndetic}
        For any $s$, $k \in \N^+$ and a finite set $\Phi \subset \G$, there is $n = n(s,k, \Phi) \in \N$ such that for all finite $F \subset \G$ with $|F| \geq n$, the following holds:

        \smallskip
        
        \noindent For every free zero-dimensional Polish $\G$-space $X$ such that $\casi(X) \leq s$, there exist disjoint $F$-syndetic clopen subsets $A_1$, \ldots, $A_k \subseteq X$ whose union 
          $A \defeq A_1 \sqcup \ldots \sqcup A_k$ is $\Phi$-spaced.
    \end{lemmacopy}
    \begin{scproof}
        Let $D \defeq \Phi^{-1}\Phi$ and $d \defeq |D|$. Fix an integer $r \in \N^+$ so large that
        \begin{equation}\label{eq:large_m}
            \mathsf{e}^{s+1} \, k \left(1 - \frac{1}{(k+1)^d}\right)^r \, r^{2(s+1)} \,\leq\, 1.
        \end{equation}
        This is possible since the left-hand side of \eqref{eq:large_m} tends to $0$ as $r \to \infty$. Now set $n \defeq d^2 r$.

        Suppose $F \subset \G$ is a finite set with $|F| \geq n$. Then $F$ has a $D$-spaced subset $R \subseteq F$ of size $|R| = r$ \cite[Lem.~4.1]{LLLDyn3}. 
        %
        Let $W \defeq DR$. Call a mapping $\phi \colon W \to (k+1)$ \emphd{acceptable} if for each $1 \leq i \leq k$, there is $\sigma \in R$ such that
        \begin{equation}\label{eq:acceptable}
            \phi(\delta \sigma) \,=\, \begin{cases}
                i &\text{if } \delta = \iden,\\
                0 &\text{if } \delta \neq \iden,
            \end{cases} \qquad \text{for all }\delta \in D.
        \end{equation}
        Let $\mathcal{P} \subseteq k^W$ be the set of all acceptable maps $\phi \colon W \to k$.

        \begin{claim}\label{claim:bound_on_prob}
            $\displaystyle 1 - \frac{|\mathcal{P}|}{(k+1)^{|W|}} \,\leq\, k \left(1 - \frac{1}{(k+1)^d}\right)^r$.
        \end{claim}
        \begin{claimproof}
            We seek an upper bound on the probability that a uniformly random map $\phi \colon W \to (k+1)$ is not acceptable. 
            The probability \eqref{eq:acceptable} holds for fixed $1 \leq i \leq k$ and a specific element $\sigma \in R$ is precisely $1/(k+1)^d$. Since the set $R$ is $D$-spaced, the conditions for different choices of $\sigma \in R$ are independent, so the probability \eqref{eq:acceptable} fails for some $1 \leq i \leq k$ and all $\sigma \in R$ is
            \[
                \left(1 - \frac{1}{(k+1)^d}\right)^r.
            \]
            Summing over all $1 \leq i \leq k$ yields the desired bound.
        \end{claimproof}

        Thanks to Claim~\ref{claim:bound_on_prob} and \eqref{eq:large_m}, we may apply Theorem~\ref{theo:cont_LLL} to obtain a continuous $\G$-equivariant map $\pi \colon X \to \Sigma(k+1, W, \mathcal{P})$. To complete the proof, we let $A_i$ be the set of all $x \in X$ such that
        \[
            (\pi(x))(\delta) \,=\, \begin{cases}
                i &\text{if } \delta = \iden,\\
                0 &\text{if } \delta \neq \iden,
            \end{cases} \qquad \text{for all }\delta \in D.
        \]
        The sets $A_i$ for $1 \leq i \leq k$ are clearly clopen and pairwise disjoint. By the definition of $\mathcal{P}$, we have $R^{-1} \cdot A_i = X$ for all $1 \leq i \leq k$, so $A_i$ is $F$-syndetic as $F \supseteq R$. Finally, if $x$, $y \in A \defeq A_1 \sqcup \ldots \sqcup A_k$ are distinct points, then $(\pi(y))(\iden) \neq 0$, and thus $y \notin (D \cdot x)$, i.e., $(\Phi \cdot x) \cap (\Phi \cdot y) = \0$, as desired.
    \end{scproof}

    \section{Finishing the proof of Theorem~\ref{theo:example}}\label{sec:final}

    Here we complete the proof that for $s \in \N^+$, $\Free(\Sep(s))$ is an \sr $\G$-space. Fix nonempty open sets $U_1$, \ldots, $U_k \subseteq ((s+1)^\G)^\N$. 
    %
    %
        As explained in \S\ref{subsubsec:combi+LLL}, Lemma~\ref{lemma:key} implies that for every large enough finite set $F \subset \G$, 
        we have 
        a continuous $\G$-equivariant map 
        \[
            \pi \colon \Free(\Sep(s)) \to \Sep(s)
        \]
        such that the sets $\pi^{-1}(U_1)$, \ldots, $\pi^{-1}(U_k)$ are $F$-syndetic. It remains to modify the map $\pi$ to make its range included in $\Free(\Sep(s))$. (This step is only needed if $\G$ has nontrivial torsion elements.) 
        
        Recall that $(\Phi_n)_{n \in \N}$ is a list of all finite subsets of $\G$ that includes each set infinitely often, and
        \[
            \Sep(s) \,=\, \prod_{n \in \N} \Sep(s, \Phi_n).
        \]
        We may assume without loss of generality that there is $N \in \N$ such that each set $U_i$ is of the form
        \[
            U_i \,=\, V_i \times \prod_{n \geq N} (s+1)^\G,
        \]
        for some nonempty open $V_i \subseteq ((s+1)^\G)^N$. By composing $\pi$ with the coordinate projection
        \[
            \Sep(s) \,=\, \prod_{n \in \N} \Sep(s,\Phi_n) \,\twoheadrightarrow\, \prod_{n < N} \Sep(s,\Phi_n),
        \]
        we obtain a continuous $\G$-equivariant map
        \[
            \rho \colon \Free(\Sep(s)) \to \prod_{n < N}  \Sep(s,\Phi_n)
        \]
        such that the sets $\rho^{-1}(V_1)$, \ldots, $\rho^{-1}(V_k)$ are $F$-syndetic. Next we observe that, since every finite subset of $\G$ appears in the sequence $(\Phi_n)_{n \in \N}$ infinitely often, there is a $\G$-equivariant homeomorphism 
        \[
            \sigma \colon \Sep(s) \to \prod_{n \geq N} \Sep(s,\Phi_n).
        \]
        The restriction of $\sigma$ to $\Free(\Sep(s))$ establishes a $\G$-equivariant homeomorphism
        \[
            \Free(\Sep(s)) \,\cong\, \Free\left(\prod_{n \geq N} \Sep(s,\Phi_n)\right).
        \]
        Now we put the maps $\rho$ and $\sigma$ together and define
        \[
            \pi' \colon \Free(\Sep(s)) \to \prod_{n < N} \Sep(s,\Phi_n) \times \Free\left(\prod_{n \geq N} \Sep(s,\Phi_n)\right)
        \]
        via $\pi'(x) \defeq (\rho(x), \sigma(x))$. It remains to note that
        \[
            \prod_{n < N} \Sep(s,\Phi_n) \times \Free\left(\prod_{n \geq N} \Sep(s,\Phi_n)\right) \,\subseteq\, \Free\left(\prod_{n \in \N} \Sep(s, \Phi_n)\right) \,=\, \Free(\Sep(s)),
        \]
        so $\pi'$ has all the desired properties.

    \section{Subshifts with minimal subdynamics}
    \label{sec:subshift}


    In this section we prove Corollary~\ref{corl:subshift}, i.e., we show that the $\G$-flow in Theorem~\ref{theo:main_prime} can be taken to be a subshift. 
    We shall use the following result of Seward and Tucker-Drob:

    \begin{theo}[{Seward--Tucker-Drob \cite{ST-D}}]\label{theo:ST-D}
        There exists a free subshift $Y \subset 2^\G$ such that every free Borel action $\G \acts X$ on a Polish space $X$ admits a Borel $\G$-equivariant map $\pi \colon X \to Y$.
    \end{theo}

    See also \cite{Ber_cont} for a shorter probabilistic proof of Theorem~\ref{theo:ST-D}. In addition to Theorem~\ref{theo:ST-D}, we need the following lemma:

    \begin{lemma}\label{lemma:transfer}
        Let $(\mathcal{F}_i)_{i \in I}$ be a collection of families of finite subsets of $\G$ and let 
        $X$ be a 
        $\G$-flow that is $\mathcal{F}_i$-minimal for all $i \in I$. If $Y$ is a $\G$-flow such that there exists a Baire-measurable $\G$-equivariant map $\pi \colon X \to Y$, then $Y$ has a subflow $Z \subseteq Y$ that is $\mathcal{F}_i$-minimal for all $i \in I$.
    \end{lemma}
    \begin{scproof}
        Call a closed set $C \subseteq Y$ \emphd{huge} if $\pi^{-1}(C)$ is comeager in $X$. Let the set of all huge closed subsets of $Y$ be denoted by $\mathcal{H}$. Note that if $C \in \mathcal{H}$, then $\gamma \cdot C \in \mathcal{H}$ as well for all $\gamma \in \G$, because
        $
            \pi^{-1}(\gamma \cdot C) = \gamma \cdot \pi^{-1}(C)
        $. Hence, the following is a closed $\G$-invariant subset of $Y$: 
        \[
            Z \,\defeq\, \bigcap_{C \in \mathcal{H}} C.
        \]
        The intersection of finitely (or even countably) many huge sets is still huge, and hence it is nonempty by the Baire category theorem applied in $X$. Since $Y$ is compact, this implies that $Z \neq \0$, so $Z$ is a subflow of $Y$. (Moreover, if $Y$ is Polish, then $Z$ itself is huge.) 
        
        
        We claim that $Z$ is $\mathcal{F}_i$-minimal for all $i \in I$. Indeed, take any $i \in I$ and suppose $U \subseteq Y$ is an open set such that $U \cap Z \neq \0$. We need to find a set $F \in \mathcal{F}_i$ with $F^{-1} \cdot U \supseteq Z$. 
        Since $Y$ is compact and Hausdorff, we may pick an open set $W \subseteq U$ such that $W \cap Z \neq \0$ and $C \defeq \overline{W} \subseteq U$. 
        The condition $W \cap Z \neq \0$ means that
        the closed set $Y \setminus W$ is not huge, i.e., $\pi^{-1}(W)$ is non-meager in $X$. As $C \supseteq W$, $\pi^{-1}(C)$ is non-meager as well. 
        The set $\pi^{-1}(C)$ is Baire-measurable, so there is a nonempty open set $V \subseteq X$ such that the symmetric difference $V \symdif \pi^{-1}(C)$ is meager.  As $X$ is $\mathcal{F}_i$-minimal, there is $F \in \mathcal{F}_i$ such that $F^{-1}\cdot V = X$, and hence $F^{-1} \cdot \pi^{-1}(C) = \pi^{-1}(F^{-1} \cdot C)$ is comeager in $X$. In other words, the closed set $F^{-1} \cdot C$ is huge, and therefore $Z \subseteq F^{-1} \cdot C \subseteq F^{-1} \cdot U$, as desired.
    \end{scproof}

    \begin{scproof}[ of Corollary~\ref{corl:subshift}]
        Let $(\mathcal{F}_n)_{n \in \N}$ be a sequence of unbounded families of finite subsets of $\G$. By Theorem~\ref{theo:main_prime}, there exists a free $\G$-flow $X$ that is $\mathcal{F}_n$-minimal for all $n \in \N$. Furthermore, the proof of Theorem~\ref{theo:main_prime} presented in \S\ref{subsubsec:ample} yields such a $\G$-flow $X$ whose underlying space is Polish. By Theorem~\ref{theo:ST-D}, there exist a free subshift $Y \subset 2^\G$ and a Borel $\G$-equivariant map $\pi \colon X \to Y$. Since $\pi$ is Borel, it is in particular Baire-measurable, so, by Lemma~\ref{lemma:transfer}, $Y$ has a subflow $Z \subseteq Y$
        that is $\mathcal{F}_n$-minimal for all $n \in \N$, as desired.
    \end{scproof}

    \section{Disjoint flows}\label{sec:disjoint}


    %



    Let $X$ be a $\G$-flow. For a point $x \in X$ and a set $U \subseteq X$, the set of \emphd{visiting times} of $x$ to $U$ is
    \[
        \Vis(x,U) \,\defeq\, \set{\gamma \in \G \,:\, \gamma \cdot x \in U}.
    \]
    Let $\mathcal{F}(U)$ denote the family of all finite sets $F \subset \G$ such that $F \subseteq \Vis(x,U)$ for some $x \in U$. The bridge between disjointness and minimal subdynamics is given by the following proposition:

    \begin{prop}\label{prop:characterize}
        A minimal $\G$-flow $Y$ is disjoint from a $\G$-flow $X$ if and only if $Y$ is $\mathcal{F}(U)$-minimal for all nonempty open sets $U \subseteq X$.
    \end{prop}
    \begin{scproof}
        Suppose first that $Y$ is $\mathcal{F}(U)$-minimal for all nonempty open sets $U \subseteq X$. Let $Z \subseteq X \times Y$ be a joining and consider arbitrary nonempty open sets $U \subseteq X$ and $V \subseteq Y$. Since $Y$ is $\mathcal{F}(U)$-minimal, there is a set $F \in \mathcal{F}(U)$ such that $V$ is $F$-syndetic. Let $x \in U$ be a point satisfying $F \subseteq \Vis(x,U)$. Since $Z$ projects onto $X$, we can find $y \in Y$ such that $(x,y) \in Z$. As $V$ is $F$-syndetic, there is $\sigma \in F$ with $\sigma \cdot y \in V$. Since $F \subseteq \Vis(x,U)$, we have $\sigma \cdot x \in U$, and hence
        \[
            \sigma \cdot (x,y) \,=\, (\sigma \cdot x, \sigma \cdot y) \,\in\, (U \times V) \cap Z.
        \]
        It follows that $Z$ is dense in $X \times Y$. Since $Z$ is closed, we conclude that $Z = X \times Y$, as desired.

        Now suppose $Y$ fails to be $\mathcal{F}(U)$-minimal for some nonempty open set $U \subseteq X$. This means that there is a nonempty open set $V \subseteq Y$ that is not $F$-syndetic for any $F \in \mathcal{F}(U)$. We will construct a joining $Z$ of $X$ and $Y$ such that $Z \cap (U \times V) = \0$, showing that $Y$ is not disjoint from $X$.

        For each point $x\in U$, we have $\Vis(x,U)^{-1} \cdot V \neq Y$, for otherwise, by the compactness of $X$, there would exist a finite set $F \subseteq \Vis(x,U)$ such that $V$ is $F$-syndetic. Therefore, we may pick a point \[y_x \,\in\, Y \setminus \big(\Vis(x,U)^{-1} \cdot V\big).\] Consider the following $\G$-invariant (but not necessarily closed) subset $A \subseteq X \times Y$:
        \[
            A \,\defeq\, \set{(\gamma \cdot x, \gamma \cdot y_x) \,:\, x \in U, \, \gamma \in \G}.
        \]
        We claim that $A \cap (U \times V) = \0$. Indeed, if $(x,y) \in A$, then there exist $x' \in U$ and $\gamma \in \G$ such that
        \[x \,=\, \gamma \cdot x' \qquad \text{and} \qquad y \,=\, \gamma \cdot y_{x'}.\] If $x \in U$, then  $\gamma \in \Vis(x', U)$, and hence $\gamma \cdot y_{x'} \notin V$, as claimed. It follows that
        \[
            Z \,\defeq\, \overline{A} \,\cup\, \big((X \setminus (\G \cdot U)) \times Y\big)
        \]
        is a subflow of $X \times Y$ avoiding $U \times V$. The first coordinate projection of $A$ includes $\G \cdot U$, which shows that $Z$ projects onto $X$. The second coordinate projection of $Z$ is a subflow of $Y$, which must be equal to $Y$ since $Y$ is minimal. Therefore, $Z$ is a joining of $X$ and $Y$, and we are done.
    \end{scproof}

    Now we proceed to prove Theorem~\ref{theo:disjoint_many}. 

    \begin{lemma}\label{lemma:recurrent}
        If $X$ is a 
        $\G$-flow with no wandering points, then for every nonempty open set $U \subseteq X$, there exists a point $x \in U$ such that the set $\Vis(x,U)$ is infinite.
    \end{lemma}
    \begin{scproof}
        Suppose not. Then we have $U = \bigcup_{F} C_F$, where the union is over all finite sets $F \subset \G$ and $C_F$ is the set of all points $x \in U$ with $\Vis(x,U) \subseteq F$. Each set $C_F$ is relatively closed in $U$, so, by the Baire category theorem, at least one of them has nonempty interior. In other words, there exist a nonempty open set $W \subseteq U$ and a finite subset $F \subset \G$ such that $\Vis(x,U) \subseteq F$ for all $x \in W$. But then every point $x \in W$ satisfies $\Vis(x,W) \subseteq \Vis(x,U) \subseteq F$, which implies that $W \cap (\gamma \cdot W) = \0$ for all $\gamma \in \G \setminus F$. Hence, $W$ is a wandering nonempty open set, which is a contradiction.
    \end{scproof}


    \begin{scproof}[ of Theorem~\ref{theo:disjoint_many}]
        Let $(X_n)_{n \in \N}$ be a countable family of Polish $\G$-flows with no wandering points. We seek a free minimal subflow of $2^\G$ that is disjoint from every $X_n$. For each $n \in \N$, we fix a countable basis $(U_{n,m})_{m \in \N}$ for the topology on $X_n$ consisting of nonempty open sets. Lemma~\ref{lemma:recurrent} implies that each family $\mathcal{F}(U_{n,m})$ is unbounded. By Corollary~\ref{corl:subshift}, there is a free subflow $Y \subset 2^\G$ that is $\mathcal{F}(U_{n,m})$-minimal for all $n$, $m \in \N$; in particular, $Y$ is minimal. Note that if $U \subseteq X_n$ is a nonempty open set, then $Y$ is $\mathcal{F}(U)$-minimal because $\mathcal{F}(U) \supseteq \mathcal{F}(U_{n,m})$ for any basic open set $U_{n,m} \subseteq U$. Therefore, Proposition~\ref{prop:characterize} implies that $X_n \perp Y$ for all $n \in \N$, as desired.
    \end{scproof}

    Next we prove Theorem~\ref{theo:disjoint_one}, restated here for ease of reference:

    \begin{theocopy}{theo:disjoint_one}
        The following statements are equivalent for a Polish $\G$-flow $X$:
        \begin{enumerate}[label=\ep{\normalfont\arabic*}]
            \item\label{item:infinite} $X$ is disjoint from some infinite $\G$-flow,
            \item\label{item:free} $X$ is disjoint from some free minimal subflow of $2^\G$,
            \item\label{item:non-wandering} $X$ has no wandering points.
        \end{enumerate}
    \end{theocopy}
    \begin{scproof}
        Implication \ref{item:free} $\Longrightarrow$ \ref{item:infinite} is obvious, while \ref{item:non-wandering} $\Longrightarrow$ \ref{item:free} is given by Theorem~\ref{theo:disjoint_many}. It remains to verify \ref{item:infinite} $\Longrightarrow$ \ref{item:non-wandering}. Assume $X$ has a wandering nonempty open subset $U \subseteq X$ and let $F \subset \G$ be a finite set such that $U \cap (\gamma \cdot U) = \0$ for all $\gamma \in \G \setminus F$. Suppose $X$ is disjoint from a $\G$-flow $Y$. We need to show that $Y$ is finite.
        
        Since $X \perp Y$, at least one of $X$ and $Y$ is minimal, and $X$ cannot be minimal because minimal $\G$-flows have no wandering points. 
        Therefore, $Y$ is minimal. By Proposition~\ref{prop:characterize}, it follows that $Y$ is $\mathcal{F}(U)$-minimal. Note that each set in $\mathcal{F}(U)$ is a subset of $F$ since $\Vis(x,U) \subseteq F$ for all $x \in U$. It follows that for any point $y \in Y$, the set $F \cdot y$ is dense in $Y$. Since finite sets in Hausdorff spaces are closed, this means that $Y = F \cdot y$, so it is finite.
    \end{scproof}

    To finish this section, we prove Proposition~\ref{prop:cosets}.

    \begin{scproof}[ of Proposition~\ref{prop:cosets}]
        Let $\Delta \leq \G$ be a subgroup and let $X$ be a $\G$-flow. We need to show that $X$ is disjoint from $(\G/\Delta)^*$ if and only if $X$ is $\Delta$-minimal. Note that $(\G/\Delta)^*$ is not minimal (the point at infinity forms a proper subflow), so if $X \perp (\G/\Delta)^*$, then $X$ must be minimal. Hence, by Proposition~\ref{prop:characterize}, $X \perp (\G/\Delta)^*$ if and only if $X$ is $\mathcal{F}(U)$-minimal for all nonempty open $U \subseteq (\G/\Delta)^*$.

        The set $\set{\Delta}$ is open and $\mathcal{F}(\set{\Delta}) = \finset{\Delta}$. Thus, if $X \perp (\G/\Delta)^*$, then $X$ is $\Delta$-minimal.

        Conversely, suppose $X$ is $\Delta$-minimal and consider any nonempty open set $U \subseteq (\G/\Delta)^*$. Pick an arbitrary point $y \in U$ and let $\Delta'$ be the stabilizer of $y$. Then $\Delta'$ is either a subgroup of $\G$ conjugate to $\Delta$ or, if $y$ is the point at infinity, $\G$ itself. In either case, $X$ is $\Delta'$-minimal (it is easy to see that a $\Delta$-minimal $\G$-flow is also minimal with respect to every conjugate of $\Delta$). Since $\Delta' \subseteq \Vis(y, U)$, we conclude that $X$ is $\mathcal{F}(U)$-minimal, as desired.
    \end{scproof}

    \section{Borel complete sections}\label{subsec:subshift_appl}

    
    
    Here we use Corollary~\ref{corl:subshift_prime} to prove Theorems~\ref{theo:DST1} and \ref{theo:DST}. We begin with a couple lemmas.

    \begin{lemma}\label{lemma:one_set}
        Let $\mathcal{F}$ be a family of finite subsets of $\G$ and let $X$ be an $\mathcal{F}$-minimal 
        $\G$-flow. If $A \subseteq X$ is a non-meager Baire-measurable set, then for some $F  \in \mathcal{F}$, $F^{-1} \cdot A$ is comeager in $X$.
    \end{lemma}
    \begin{scproof}
        Let $U \subseteq X$ be the nonempty open set such that $U \symdif A$ is meager. By the $\mathcal{F}$-minimality of $X$, there is $F \in \mathcal{F}$ such that $F^{-1} \cdot U = X$, and hence $F^{-1} \cdot A$ is comeager, as desired.
    \end{scproof}

    \begin{lemma}\label{lemma:many_sets}
        Let $\mathcal{F} = \set{F_n}_{n \in \N}$ be a family of finite subsets of $\G$ and let $X$ be an $\mathcal{F}$-minimal 
        $\G$-flow. If $(A_n)_{n \in \N}$ is a sequence of non-meager Baire-measurable subsets of $X$, then the set $\bigcup_{n \in \N} (F_n \cdot A_n)$ is comeager in $X$.
    \end{lemma}
    \begin{scproof}
        Suppose for contradiction that the set $M \defeq X \setminus \left(\bigcup_{n \in \N} (F_n \cdot A_n)\right)$ is non-meager. Then, by Lemma~\ref{lemma:one_set}, there is $n \in \N$ such that $F_n^{-1} \cdot M$ is comeager in $X$. This implies that $(F_n^{-1} \cdot M) \cap A_n \neq \0$, i.e., $M \cap (F_n \cdot A_n) \neq \0$, which contradicts the definition of $M$. 
    \end{scproof}

    We are now ready to deduce Theorems~\ref{theo:DST1} and \ref{theo:DST}, restated here for ease of reference.

    \begin{theocopy}{theo:DST1}\label{theo:DST1_copy}
        Suppose $(B_n)_{n \in \N}$ is a sequence of Borel complete sections in $\Free(2^\G)$ and $(F_n)_{n \in \N}$ is a sequence of finite subsets of $\G$ such that $\sup_{n \in \N}|F_n| = \infty$. Then there exists a point $x \in \Free(2^\G)$ satisfying $x \in F_n \cdot B_n$ for infinitely many $n \in \N$.
    \end{theocopy}
    \begin{scproof}
        For each $m \in \N$, let $\mathcal{F}_m \defeq \set{F_n}_{n \geq m}$. Since the families $\mathcal{F}_m$ are unbounded, Corollary~\ref{corl:subshift_prime} yields a free subshift $X \subset \Free(2^\G)$ that is $\mathcal{F}_m$-minimal for all $m \in \N$. Let $A_n \defeq B_n \cap X$. Since $B_n$ is a complete section, we have $\G \cdot A_n = X$, so $A_n$ is non-meager in $X$. By Lemma~\ref{lemma:many_sets}, we conclude that for each $m \in \N$, the set $\bigcup_{n \geq m} (F_n \cdot A_n)$ is comeager in $X$. Thus, the set $\bigcap_{m \in \N} \bigcup_{n \geq m} (F_n \cdot A_n)$ of all points $x \in X$ satisfying $x \in F_n \cdot B_n$ for infinitely many $n \in \N$ is also comeager in $X$. In particular, this set is nonempty. 
    \end{scproof}

    \begin{theocopy}{theo:DST}
        If $B \subseteq \Free(2^\G)$ is a Borel complete section, then there exists $n \in \N$ such that for every finite set $F \subset \G$ of size at least $n$, the set $F \cdot B$ traps a point.
    \end{theocopy}
    \begin{scproof}
        Suppose for contradiction that for each $n \in \N$, there exists a set $F_n \subset \G$ of size $|F_n| = n$ such that $F_n \cdot B$ does not trap a point. Let $\mathcal{F} \defeq \set{F_n^{-1}}_{n \in \N}$ and consider a free $\mathcal{F}$-minimal subshift $X \subset \Free(2^\G)$. As $B$ is a complete section, the set $A \defeq B \cap X$ is non-meager in $X$, so, by Lemma~\ref{lemma:one_set}, $F_n \cdot A$ is comeager in $X$ for some $n \in \N$. Since $F_n \cdot A$ does not trap a point by assumption, the set $M \defeq X \setminus (F_n \cdot A)$ must meet every orbit in $X$, which is impossible as $M$ is meager in $X$.
    \end{scproof}

    Note that our proof of Theorem~\ref{theo:DST1_copy} actually yields the following somewhat stronger statement: For any family $(B_{n,m})_{n,m \in \N}$ of Borel complete sections in $\Free(2^\G)$ and any collection $(F_{n,m})_{n,m \in \N}$ of finite subsets of $\G$ satisfying $\sup_{n \in \N} |F_{n,m}| = \infty$ for all $m \in \N$, the following set traps a point:
    \[
        \bigcap_{m \in \N} \bigcup_{n \in \N} (F_{n,m} \cdot B_{n,m}).
    \]

    \subsection*{Acknowledgments}

    We are grateful to Be'eri Greenfeld, Gil Goffer, and Tianyi Zheng for organizing the AMS Special Session on \emph{Combinatorial and Probabilistic Methods in Group Theory} at the 2025 Joint Mathematics Meetings (January 8--11, 2025) in Seattle, where this project began, and to Lewis Bowen for providing a stimulating environment during the UT Austin Mini-school in \emph{Analytic Group Theory} (March 17--21, 2025), where part of our work on this project took place. We also thank Andy Zucker for insightful conversations. This material is based upon work partially supported by the Alfred P.~Sloan Foundation and the National Science Foundation under grants DMS-2528522 and DMS-2348981. Any opinions, findings, and conclusions or recommendations expressed in this material are those of the authors and do not necessarily reflect the views of the National Science Foundation.

\printbibliography

@article{FSZ,
	author = {J. Frisch and B. Seward and A. Zucker},
	title = {Minimal subdynamics and minimal flows without characteristic measures},
	journaltitle = {Forum Math. Sigma},
	volume = {12},
	pages = {\#{}e58},
	date = {2024},
}

@article{Zucker,
	author = {A. Zucker},
	title = {Minimal flows with arbitrary centralizer},
	journaltitle = {Ergod. Theory Dyn. Syst.},
	volume = {42},
    number = {1},
	pages = {310--320},
	date = {2022},
}

@article{CJMSTD,
  title={Borel asymptotic dimension and hyperfinite equivalence relations},
  author = {Conley, C.T. and Jackson, S. and Marks, A.S. and Seward, B. and Tucker-Drob, R.D.},
  journal={Duke Math. J.},
  volume={172},
  number={16},
  pages={3175--3226},
  year={2023},
}

@article{EL,
	author = {P. Erd{\H{o}}s and L. Lov{\'{a}}sz},
	title = {Problems and results on $3$-chromatic hypergraphs and some related questions},
	journaltitle = {Infinite and Finite Sets, Colloq. Math. Soc. J. Bolyai},
	editor = {A. Hajnal and R. Rado  and V.T. S{\'{o}}s},
	publisher = {North Holland},
	date = {1975},
	pages = {609--627},
}

@article{BW,
	author = {A. Bernshteyn and F. Weilacher},
	title = {Borel versions of the Local Lemma and $\mathsf{LOCAL}$ algorithms for graphs of finite asymptotic separation index},
    journaltitle = {Trans. Am. Math. Soc.},
	date = {2025},
    addendum = {Published electronically, preprint at \url{https://arxiv.org/abs/2308.14941}},
    url = {https://doi.org/10.1090/tran/9455},
}

@article{FTVF,
	author = {J. Frisch and O. Tamuz and P. Vahidi Ferdowsi},
	title = {Strong amenability and the infinite conjugacy class property},
	journaltitle = {Invent. Math.},
    volume = {218},
	date = {2019},
	pages = {833--851},
}

@article{FT,
	author = {J. Frisch and O. Tamuz},
	title = {Symbolic dynamics on amenable groups: the entropy of generic shifts},
	journaltitle = {Ergod. Theory Dyn. Syst.},
    volume = {37},
    number = {4},
	date = {2017},
	pages = {1187--1210},
}

@book{subshifts1,
    author = {D. Lind and B. Marcus},
    title = {An Introduction to Symbolic Dynamics and Coding},
    publisher = {Cambridge University Press},
    date = {1995},
}

@book{subshifts2,
    author = {T. Ceccherini-Silberstein and M. Coornaert},
    title = {Cellular Automata and Groups},
    location = {Berlin, Heidelberg},
    publisher = {Springer},
    date = {2010},
}

@article{KST,
  title={Borel chromatic numbers},
  author={A.S. Kechris and S. Solecki and S. Todorcevic},
  journal={Adv. Math.},
  volume={141},
  number={1},
  pages={1--44},
  date={1999},
}

@article{Ber_cont,
	author = {A. Bernshteyn},
	title = {Probabilistic constructions in continuous combinatorics and a bridge to distributed algorithms},
	date = {2023},
    journaltitle = {Adv. Math.},
    volume = {415},
    pages = {\#{}108895},
}

@unpublished{KechrisMarks,
	author = {A.S. Kechris and A.S. Marks},
	title = {Descriptive Graph Combinatorics},
	date = {2020},
	howpublished = {\url{https://math.berkeley.edu/~marks/papers/combinatorics20book.pdf} (preprint)},
}

@article{Pikh_survey,
	author = {O. Pikhurko},
	title = {Borel combinatorics of locally finite graphs},
	journaltitle = {Surveys in Combinatorics, 28th British Combinatorial Conference},
	editor = {K.K. Dabrowski {et al.}},
	pages = {267--319},
	date = {2021},
}

@article{Notices,
    author = {A. Bernshteyn},
	title = {Descriptive combinatorics and distributed algorithms},
	date = {2022},
	journaltitle = {Not. Am. Math. Soc.},
	volume = {69},
    number = {9},
	pages = {1496--1507},
}

@book{KechrisDST,
	author = {A.S. Kechris},
	title = {Classical Descriptive Set Theory},
	date = {1995},
	publisher = {Springer-Verlag},
	location = {New York},
}

@book{gromov1993asymptotic,
	author = {M. Gromov},
	title = {Geometric group theory, Vol. 2: Asymptotic invariants of infinite groups},
	date = {1993},
	publisher = {Cambridge Univ. Press},
	location = {Cambridge},
	editor = {A. Niblo and M.A. Roller},
}

@unpublished{ASIalgorithms,
	author = {L. Qian and F. Weilacher},
	title = {Descriptive combinatorics, computable combinatorics, and ASI algorithms},
	howpublished = {\url{https://arxiv.org/abs/2206.08426} (preprint)},
	date = {2022},
}

@article{BWKonig,
	author = {M. Bowen and F. Weilacher},
	title = {Definable K{\H{o}}nig Theorems},
	journaltitle = {Proc. Am. Math. Soc.},
    volume = {151},
    date = {2023},
    pages = {4991--4996},
}

@article{WeilacherFinDim,
  title={Borel edge colorings for finite-dimensional groups},
  author={F. Weilacher},
  journal={Isr. J. Math.},
  pages={737--780},
  volume = {263},
  date={2024},
}

@unpublished{Grids2,
	author = {A. Bernshteyn and J. Yu},
	title = {Embedding Borel graphs into grids of asymptotically optimal dimension},
	howpublished = {\url{https://arxiv.org/abs/2407.19785} (preprint)},
	date = {2024},
}

@book{AS,
	author = {N. Alon and J.H. Spencer},
	title = {The Probabilistic Method},
	date = {2016},
	edition = {4},
	publisher = {John Wiley {\&} Sons},
}

@book{MolloyReed,
	author = {M. Molloy and B. Reed},
	title = {Graph Colouring and the Probabilistic Method},
	publisher = {Springer-Verlag},
	location = {Berlin Heidelberg},
    date = {2002},
}

@article{ABT,
  title={Realization of aperiodic subshifts and uniform densities in groups},
  author={N. Aubrun and S. Barbieri and S. Thomass{\'{e}}},
  journal={Groups Geom. Dyn.},
  pages={107--129},
  volume = {13},
number = {1},
  date={2019},
}

@article{Elek,
  title={Uniformly recurrent subgroups and simple $C^*$-algebras},
  author={G. Elek},
  journal={J. Funct. Anal.},
  pages={1657--1689},
  volume = {274},
number = {6},
  date={2018},
}

@article{LLLDyn1,
    author = {A. Bernshteyn},
	title = {Building large free subshifts using the Local Lemma},
	date = {2019},
	journaltitle = {Groups Geom. Dyn.},
	volume = {13},
    number = {4},
	pages = {1417--1436},
}

@article{LLLDyn2,
    author = {A. Bernshteyn},
	title = {Ergodic theorems for the shift action and pointwise versions of the Ab{\'{e}}rt--Weiss theorem},
	date = {2020},
	journaltitle = {Isr. J. Math.},
	volume = {235},
	pages = {255--293},
}

@article{LLLDyn3,
    author = {A. Bernshteyn},
	title = {A short proof of Bernoulli disjointness via the Local Lemma},
	date = {2020},
	journaltitle = {Proc. Am. Math. Soc.},
	volume = {148},
    number = {12},
	pages = {5235--5240},
}

@article{LLLDyn4,
    author = {A. Bernshteyn},
	title = {Equivariant maps to subshifts whose points have small stabilizers},
	date = {2023},
	journaltitle = {J. Mod. Dyn.},
	volume = {19},
	pages = {1--30},
}

@article{SpencerRamsey,
	author = {J. Spencer},
	title = {Asymptotic lower bounds for Ramsey functions},
	date = {1977},
    journaltitle = {Disc. Math.},
    volume = {20},
    pages = {69--76},
}

@article{Beck,
	author = {J. Beck},
	title = {An algorithmic approach to the Lov{\'{a}}sz Local Lemma. I},
	journaltitle = {Rand. Str. {\&} Alg.},
	volume = {2},
	number = {4},
	pages = {343--365},
	date = {1991},
}

@article{BGRDeterministicLLL,
	author = {S. Brandt and C. Grunau and V. Rozho{\v{n}}},
	title = {Generalizing the sharp threshold phenomenon for the distributed complexity of the Lov{\'{a}}sz Local Lemma},
	journaltitle = {ACM Symposium on Principles of Distributed Computing (PODC)},
	date = {2020},
	pages = {329--338},
	addendum = {Full version: \url{https://arxiv.org/abs/2006.04625}},
}

@article{BMUDeterministicLLL,
	author = {S. Brandt and Y. Maus and J. Uitto},
	title = {A sharp threshold phenomenon for the distributed complexity of the Lov{\'{a}}sz Local Lemma},
	journaltitle = {ACM Symposium on Principles of Distributed Computing (PODC)},
	date = {2019},
	pages = {389--398},
	addendum = {Full version: \url{https://arxiv.org/abs/1908.06270}},
}

@article{FG,
	author =	{M. Fischer and M. Ghaffari},
	title =	{{Sublogarithmic distributed algorithms for Lov{\'a}sz Local Lemma, and the complexity hierarchy}},
	journaltitle =	{International Symposium on DIStributed Computing (DISC)},
	pages =	{18:1--18:16},
	year =	{2017},
	volume =	{91},
	addendum = {Full version: \url{https://arxiv.org/abs/1705.04840}},
}

@article{MoserTardos,
	author = {R. Moser and G. Tardos},
	title = {A constructive proof of the general Lov\'{a}sz Local Lemma},
	journaltitle = {J. ACM},
	date = {2010},
	volume = {57},
	number = {2},
}

@article{RSh,
	author = {A. Rumyantsev and A. Shen},
	title = {Probabilistic constructions of computable objects and a computable version of Lov{\'{a}}sz Local Lemma},
	journaltitle = {Fundamenta Informaticae},
	volume = {132},
	date = {2014},
	number = {1},
	pages = {1--14},
}

@unpublished{CGMPT,
	author = {E. Cs{\'{o}}ka and \L{}. Grabowski and A. M{\'{a}}th{\'{e}} and O. Pikhurko and K. Tyros},
	title = {Moser--Tardos Algorithm with small number of random bits},
	date = {2024},
	howpublished = {\url{https://arxiv.org/abs/2203.05888} (preprint)},
}

@article{BerDist,
	author = {A. Bernshteyn},
	title = {Distributed algorithms, the Lov{\'{a}}sz Local Lemma, and descriptive combinatorics},
	date = {2023},
	journaltitle = {Invent. Math.},
	volume = {233},
	pages = {495--542},
}

@article{BerMeas,
	author = {A. Bernshteyn},
	title = {Measurable versions of the Lov{\'{a}}sz Local Lemma and measurable graph colorings},
	date = {2019},
	journaltitle = {Adv. Math.},
	volume = {353},
	pages = {153--223},
}

@unpublished{BerYuBorelLLL,
	author = {A. Bernshteyn and J. Yu},
	title = {Borel Local Lemma: arbitrary random variables and limited exponential growth},
	date = {2024},
	howpublished = {\url{https://arxiv.org/abs/2412.11571} (preprint)},
}

@book{RandAlg,
	author = {R. Motwani and P. Raghavan},
	title = {Randomized Algorithms},
	date = {1995},
	publisher = {Cambridge University Press},
}

@article{ST-D,
	author = {B. Seward and R.D. Tucker-Drob},
	title = {Borel structurability on the 2-shift of a countable group},
	journaltitle = {Ann. Pure Appl. Log.},
	volume = {167},
    number = {1},
	pages = {1--21},
	date = {2016},
}

@article{GJS1,
	author = {S. Gao and S. Jackson and B. Seward},
	title = {A coloring property for countable groups},
	journaltitle = {Math. Proc. Camb. Philos. Soc.},
	volume = {147},
	date = {2009},
	pages = {579--592},
}

@article{GJS2,
	author = {S. Gao and S. Jackson and B. Seward},
	title = {Group colorings and Bernoulli subflows},
	journaltitle = {Mem. Am. Math. Soc.},
	volume = {241},
	number = {1141},
	date = {2016},
}

@article{DS,
	author = {A. Dranishnikov and V. Schroeder},
	title = {Aperiodic colorings and tilings of Coxeter groups},
	journaltitle = {Groups Geom. Dyn.},
	volume = {1},
	date = {2007},
	pages = {301--318},
}

@article{GU,
	author = {E. Glasner and V.V. Uspenskij},
	title = {Effective minimal subflows of Bernoulli flows},
	journaltitle = {Proc. Am. Math. Soc.},
	volume = {137},
	date = {2009},
	pages = {3147--3154},
}

@article{Marks,
	author = {A.S. Marks},
	title = {A determinacy approach to Borel combinatorics},
	journaltitle = {J. Am. Math. Soc.},
	volume = {29},
	date = {2016},
	pages = {579--600},
}

@article{HedlundSturmian,
	author = {G.A. Hedlund},
	title = {Sturmian minimal sets},
	journaltitle = {Am. J. Math.},
	volume = {66},
    number = {4},
	date = {1944},
	pages = {605--620},
}

@book{GottHed,
	author = {W.H. Gottschalk and G.A. Hedlund},
	title = {Topological Dynamics},
	publisher = {Am. Math. Soc.},
	location = {Providence},
	date = {1955},
}

@book{RSErgodic,
	author = {S. Rubinstein-Salzedo},
	title = {Ergodic Theory},
	publisher = {Am. Math. Soc.},
	location = {Providence},
	date = {2025},
}

@article{Expansive1,
	author = {M. Boyle and D. Lind},
	title = {Expansive subdynamics},
	journaltitle = {Trans. Am. Math. Soc.},
	volume = {349},
    number = {1},
	date = {1997},
	pages = {55--102},
}

@article{Expansive2,
	author = {M. Einsiedler and D. Lind and R. Miles and T. Ward},
	title = {Expansive subdynamics for algebraic $\Z^d$-actions},
	journaltitle = {Ergod. Theory Dyn. Syst.},
	volume = {21},
    number = {6},
	date = {2001},
	pages = {1695--1729},
}

@article{Expansive3,
	author = {V. Cyr and D. Kra},
	title = {Nonexpansive $\Z^2$-subdynamics and Nivat's conjecture},
	journaltitle = {Trans. Am. Math. Soc.},
	volume = {367},
    number = {9},
	date = {2015},
	pages = {6487--6537},
}

@article{Expansive4,
	author = {M. Hochman},
	title = {On the dynamics and recursive properties of multidimensional symbolic systems},
	journaltitle = {Invent. Math.},
	volume = {176},
	date = {2009},
	pages = {131--167},
}

@article{Expansive5,
	author = {R. Pavlov and M. Schraudner},
	title = {Classification of sofic projective subdynamics of multidimensional shifts of finite type},
	journaltitle = {Trans. Am. Math. Soc.},
	volume = {367},
    number = {5},
	date = {2015},
	pages = {3371--3421},
}

@article{Expansive6,
	author = {C.F. Colle},
	title = {On periodic decompositions, one-sided nonexpansive directions and Nivat's conjecture},
	journaltitle = {Discrete Contin. Dyn. Syst.},
	volume = {43},
    number = {12},
	date = {2023},
	pages = {4299--4327},
}

@article{Expansive7,
	author = {V. Salo},
	title = {Subshifts with sparse traces},
	journaltitle = {Studia Math.},
	volume = {255},
    number = {2},
	date = {2020},
	pages = {159--207},
}

@article{Expansive8,
	author = {M.H. Schraudner},
	title = {One-dimensional projective subdynamics of uniformly mixing $\Z^d$ shifts of finite type},
	journaltitle = {Ergod. Theory Dyn. Syst.},
	volume = {35},
    number = {6},
	date = {2015},
	pages = {1962--1999},
}

@article{Expansive9,
	author = {M. Hochman},
	title = {Non-expansive directions for $\Z^2$ actions},
	journaltitle = {Ergod. Theory Dyn. Syst.},
	volume = {31},
    number = {1},
	date = {2011},
	pages = {91--112},
}

@article{GJKS_continuous,
	author = {S. Gao and S. Jackson and E. Krohne and B. Seward},
	title = {Continuous combinatorics of abelian group actions},
	journaltitle = {Mem. Am. Math. Soc.},
	volume = {311},
    number = {1573},
	date = {2025},
}

@article{GJKS_Borel,
	author = {S. Gao and S. Jackson and E. Krohne and B. Seward},
	title = {Forcing constructions and countable Borel equivalence relations},
	journaltitle = {J. Symb. Log.},
	volume = {87},
    number = {3},
	date = {2022},
    pages = {873--893},
}

@article{BernoulliDisjointness,
	author = {E. Glasner and T. Tsankov and B. Weiss and A. Zucker},
	title = {Bernoulli disjointness},
	journaltitle = {Duke Math. J.},
	volume = {170},
    number = {4},
	date = {2021},
    pages = {615--651},
}

@article{OxtobyUlam,
	author = {J.C. Oxtoby and S.M. Ulam},
	title = {Measure-preserving homeomorphisms and metrical transitivity},
	journaltitle = {Ann. Math.},
	volume = {42},
    number = {4},
	date = {1941},
    pages = {874--920},
}

@article{Halmos1,
	author = {P.R. Halmos},
	title = {Approximation theories for measure preserving transformations},
	journaltitle = {Trans. Am. Math. Soc.},
	volume = {55},
    number = {1},
	date = {1944},
    pages = {1--18},
}

@article{Halmos2,
	author = {P.R. Halmos},
	title = {In general a measure preserving transformation is mixing},
	journaltitle = {Ann. Math.},
	volume = {45},
    number = {4},
	date = {1944},
    pages = {786--792},
}

@article{Rokhlin,
	author = {V.A. Rokhlin},
    title = {{\cyrins{Общее преобразование с инвариантной мерой не есть перемешивание}} \normalfont{({Russian}) [A general measure preserving transformation is not mixing]}},
	journaltitle = {Proc. USSR Acad. Sci.},
	volume = {60},
    number = {3},
	date = {1948},
    pages = {349--351},
}

@article{Gen1,
	author = {A. Artigue},
	title = {Generic dynamics on compact metric spaces},
	journaltitle = {Topol. Its Appl.},
	volume = {255},
	date = {2019},
    pages = {1--14},
}

@article{Gen2,
	author = {S. Iyer and F. Shinko},
	title = {Asymptotic dimension and hyperfiniteness of generic Cantor actions},
	journaltitle = {Groups Geom. Dyn.},
	date = {2025},
    addendum = {Published online first, preprint at \url{https://arxiv.org/abs/2409.03078}},
}

@unpublished{Gen3,
	author = {J. Frisch and A. Kechris and F. Shinko and Z. Vidny{\'{a}}nszky},
	title = {Realizations of countable Borel equivalence relations},
	date = {2025},
	howpublished = {\url{https://arxiv.org/abs/2109.12486} (preprint)},
}

@article{Gen4,
	author = {N.C. Bernardes Jr. and U.B. Darji},
	title = {Graph theoretic structure of maps of the Cantor space},
	journaltitle = {Adv. Math.},
	volume = {231},
    number = {3--4},
	date = {2012},
    pages = {1655--1680},
}

@article{Gen5,
	author = {R. Pavlov and S. Schmieding},
	title = {On the structure of generic subshifts},
	journaltitle = {Nonlinearity},
	volume = {36},
    number = {9},
	date = {2023},
    pages = {4904--4953},
}

@article{Gen6,
	author = {M. Hochman},
	title = {Genericity in topological dynamics},
	journaltitle = {Ergod. Theory Dyn. Syst.},
	volume = {28},
    number = {1},
	date = {2008},
    pages = {125--165},
}

@article{Gen7,
	author = {E. Akin and E. Glasner and B. Weiss},
	title = {Generically there is but one self homeomorphism of the Cantor set},
	journaltitle = {Trans. Am. Math. Soc.},
	volume = {360},
	date = {2008},
    pages = {3613--3630},
}

@article{Gen8,
	author = {E. Akin and M. Hurley and J.A. Kennedy},
	title = {Dynamics of topologically generic homeomorphisms},
	journaltitle = {Mem. Am. Math. Soc.},
	volume = {167},
    number = {783},
	date = {2003},
}

@article{Gen9,
	author = {A.S. Kechris and C. Rosendal},
	title = {Turbulence, amalgamation, and generic automorphisms of homogeneous structures},
	journaltitle = {Proc. Lond. Math. Soc.},
	volume = {94},
    number = {2},
	date = {2006},
    pages = {302--350},
}

@article{Thue,
	author = {A. Thue},
	title = {{\"{U}}ber unendliche Zeichenreihen \normalfont{({German}) [On infinite strings of symbols]}},
	journaltitle = {Norske Vid. Selsk. Skr. I. Mat. Nat. Kl. Christiania},
	volume = {7},
	date = {1906},
    pages = {1--22},
}

@article{Morse,
	author = {H.M. Morse},
	title = {Recurrent geodesics on a surface of negative curvature},
	journaltitle = {Trans. Am. Math. Soc.},
	volume = {22},
    number = {1},
	date = {1921},
    pages = {84--100},
}

@article{MorseHedlund,
	author = {M. Morse and G.A. Hedlund},
	title = {Unending chess, symbolic dynamics and a problem in semigroups},
	journaltitle = {Duke Math. J.},
	volume = {11},
	date = {1944},
    pages = {1--7},
}

@article{Furstenberg,
	author = {H. Furstenberg},
	title = {Disjointness in ergodic theory, minimal sets, and a problem in Diophantine approximation},
	journaltitle = {Math. Syst. Theory},
	volume = {1},
	date = {1967},
    pages = {1--49},
}

@article{Rosenthal,
	author = {A. Rosenthal},
	title = {On strictly ergodic models for commuting ergodic transformations},
	journaltitle = {Ann. Henri Poincare},
	volume = {25},
    number = {1},
	date = {1989},
    pages = {73--92},
}

\end{document}